\documentclass[12pt]{article}
\newtheorem{thm}{Theorem}[section]

\newtheorem{prop}[thm]{Proposition}
\newtheorem{cor}[thm]{Corollary}

\newtheorem{ques}[thm]{Question}

\newcommand{\R}{\hbox{\bf R}}

\newcommand{\beeq}{\begin{eqnarray*}}
\newcommand{\eneq}{\end{eqnarray*}}
\newcommand{\proof}{\noindent {\it Proof.\hspace{4mm}}}

\newcommand{\qfd}{\hfill $\fbox{}$\vspace{4mm}}\def\newpic#1{%
\def\emline##1##2##3##4##5##6{%
\put(##1,##2){\special{em:point #1##3}}%
\put(##4,##5){\special{em:point #1##6}}%
\special{em:line #1##3,#1##6}}}
\newpic{}
\def\emline#1#2#3#4#5#6{%
\put(#1,#2){\special{em:moveto}}%
\put(#4,#5){\special{em:lineto}}}
\def\newpic#1{}
\newcommand\ZZ{{{\rm Z}\kern-.28em{\rm Z}}}

\title{Perfect domination in rectangular grid graphs}
\author{Italo J. Dejter
\\ University of Puerto Rico \\ Rio Piedras, PR 00931-3355 \\ ijdejter@uprrp.edu
\and Abel A. Delgado
\\ Auburn University \\ Auburn, AL 36849-5310 \\ delgaaa@auburn.edu}\date{}

\begin{document}
\maketitle

\begin{abstract}
A dominating set $S$ in a graph $G$ is said to be perfect if every
vertex of $G$ not in $S$ is adjacent to just one vertex of $S$.
Given a vertex subset $S'$ of a side $P_m$
of an $m\times n$ grid graph $G$, the perfect dominating sets $S$ in $G$ with $S'=S\cap V(P_m)$ can be determined via an exhaustive
algorithm $\Theta$ of running time $O(2^{m+n})$. Extending
$\Theta$ to infinite grid graphs of width $m-1$, periodicity makes
the binary decision tree of $\Theta$ prunable into a finite threaded
tree, a closed walk of which yields all such sets $S$.
The graphs induced by the complements of such sets $S$
can be codified by arrays of
ordered pairs of positive integers via $\Theta$, for the growth and
determination of which a speedier %greedy
algorithm exists. %and their periodic structure, further studied.
A recent characterization of grid graphs having total perfect
codes $S$ (with just 1-cubes as induced components), due to
Klostermeyer and Goldwasser, is given in terms of $\Theta$, which
allows to show that these sets $S$ are restrictions of only one
total perfect code $S_1$ in the integer lattice graph ${\Lambda}$ of $\R^2$.
Moreover, the complement ${\Lambda}-S_1$ yields an aperiodic
tiling, like the Penrose tiling. In contrast, the parallel,
horizontal, total perfect codes in ${\Lambda}$ are in 1-1
correspondence with the doubly infinite $\{0,1\}$-sequences.
\end{abstract}

\section{Introduction}

The {\it integer lattice graph} ${\Lambda}$ of $\R^2$ is the graph with vertex
set $\{(i,j):i,j\in\ZZ\}$ and such that any two vertices of
${\Lambda}$ are adjacent if and only if their Euclidean
distance is 1.
${\Lambda}$ and its subgraphs are represented orthogonally, their
vertical paths from left
to right for increasing indices $i\in[0,m)$ and its horizontal paths
downward for increasing indices $j=0,\dots,n-1$.

A vertex subset $S$ in a graph $G$ is said to be a
{\it perfect dominating set} (PDS) in $G$ if each vertex of
the complementary graph $G\setminus S$ of $S$ in $G$
is adjacent to just one element of $S$, (see \cite{W}). The
NP-completeness of finding an isolated PDS in $G$ as well as that
of finding a minimal PDS in a planar graph were established
respectively in \cite{BBS,K} and in Sections 3 and 4 of \cite{FH},
even if its induced components are $i$-cubes with $i\leq 1$. Thus,
it makes sense to consider the problem of existence of a PDS $S$
in a graph $G$ under an initial condition given by a vertex subset
$S'\subseteq S$ of a fixed subgraph $G'$ of $G$, that is $S\cap
G'=S\cap V(G')=S'$.

Motivated by Theorem 1 of \cite{DD1}, that says that there is no
algorithmic characterization of PDSs in ${\Lambda}$,
we consider the problem above for the case
in which $G$ is a (rectangular) grid graph with $G'$
equal to a side, (maximum lateral path), of $G$. For this case, we present
an algorithm as claimed in the Abstract.
In the rest of Section 1, we present
preliminary concepts and results and an overview of the remaining sections.

\subsection{PDSs with initial conditions in a (periodic) graph}

\begin{prop}
The following conditions are equivalent, as well as necessary, for
the existence of a {\rm PDS} $S$ in $G$ such that $S\cap G'=S'$:
\begin{enumerate}
\item[\rm 1.] No two components of the induced subgraph $G'[S']$
are at distance $2$.
\item[\rm 2.] $G'[S']$ has its components separated by distances $\geq 3$.
\item[\rm 3.] There is a unique subgraph $G''$ of $G'$ such that $S'$ is a
{\rm PDS} in $G''$ satisfying $S\cap G''=S'$.
\end{enumerate}\end{prop}

\proof Assuming that there are two components $C_1, C_2$ of
$G'[S']$ at distance 2, the unique vertex of $G'$ internal in a
path realizing that distance is dominated by respective boundary
vertices in $C_1,C_2$, which should be vertices of any PDS $S$ in $G$
such that $S\cap G'=S'$, contrary to the definition of a PDS. This
yields item 1, which in turn is seen to be equivalent to items 2
and 3. In fact, if $S$ is a PDS in $G$ such that $S\cap G'=S'$,
then let $S''$ be the set of vertices of $G'$ either dominated by
or belonging to $S'$. Then clearly $G''=G'[S'']$ is as in item 3. \qfd

A vertex subset $S'$ of $G'$ satisfying any of the three
conditions in Proposition 1.1 is said to be {\it admissible}. We
deal with the following questions.

\begin{ques}
Given connected graphs $G'\subseteq G$ and an admissible vertex
subset $S'$ of $G'$, does there exist a {\rm PDS} $S$ in $G$ such
that $S\cap G'=S'$?
\end{ques}

\begin{ques} Let $G$ be the union of connected graphs $G'=G_1\subset
G_2\subset\ldots\subset G_n\subset\ldots$. Let $S'$ be an
admissible vertex subset of $G'$. Does there exist a {\rm PDS} $S$
of some $G_n$ such that $S\cap G'=S'$?
\end{ques}

A particular case of Question 1.3 is posed in Question 4.1 as a result
of the grid-graph setting of Subsection 1.2 below.

The following fact will be used from Section 4 on. A graph $G$ as in
Question 1.3 is said to be {\it periodic} if there exists a graph
isomorphism $\eta_i:H_i=G_{i+1}-V(G_i)\rightarrow G_1$ extendible to
a graph isomorphism $G_{i+2}-V(G_i)\rightarrow G_2$, for each $i>0$,
ensuring that the attachment of $H_i$ to $G_i$ does not depend on
$i$. A PDS $S$ in a periodic graph $G$ is said to be {\it periodic}
if there are positive integers $k$ and $\ell$ such that
$\eta_{i+\ell}^{-1}\circ\eta_i(S\cap H_i)=S\cap H_{i+\ell}$ , for
each $i\ge k$, so that the PDS-{\it slices} $S\cap H_i$ and $S\cap
H_{i+\ell}$ are equivalently disposed in $H_i$ and $H_{i+\ell}$, respectively.

\begin{thm}
Given a periodic graph $G$ with $G'$ finite and an admissible
vertex subset $S'$ of $G'$, there exist a periodic {\rm PDS} $S$ in $G$
such that $S\cap G'=S'$.
\end{thm}

\proof Since $G_1=G'$ is finite, there is just a finite number of
candidates for $\eta_i(S\cap H_i)\subseteq G'$, where $0<i\in\ZZ$.
We construct a PDS $S$ in $G$ such that $S\cap G'=S'$, where $S$ is the union
of subsets $S_j$ of $G$ such that $S_j\cap G'=S'$, $S_j\subset
S_{j+1}$ and $S_{j+1}\setminus S_j=\{v_j\}$, for some vertex $v_j$
of $G$, where $j=0,1,2,\ldots$ and $S_0=S'$. Such vertices $v_j$
are referred to as $S$-{\it candidates}, some of which are eligible to
form part of a PDS $S$ in $G$, while some others are compulsorily its
members. Any assumed linear order for the vertices of $G'$
induces a linear order for the vertices of $H_i$ via $\eta_i$.
By means of such an order, select
successively vertices $v_j$ in $H_1$, continuing in $H_2$, etc., in the greediest fashion.
In the limit, yielding
the union of all $S_j$, the sought PDS $S$ is attained.
Because of the finiteness of the number or candidates for $\eta_i(S\cap
H_i)\subseteq G'$, there must exist two positive integers $x,y$
such that $x<y$ and $\eta_x(S\cap H_x)=\eta_y(S\cap H_y)$. Since
each $\eta_i$ is extendible to a graph isomorphism
$G_{i+1}-V(G_{i-1})\rightarrow G_2$, then $\eta_{x+1}(S\cap
H_{x+1})=\eta_{y+1}(S\cap H_{y+1})$. As a result, there are two
positive integer $k,\ell$ as in the definition of a periodic PDS.
For example, $k=x$ and $\ell=y-x$. \qfd

\subsection{PDSs in grid graphs}

Let $P_r$ denote a path of length $r-1$, ($2\leq
r\in\ZZ)$. The notation $P_r$ extends to $P_1=\{v\}$, formed by
a vertex $v$, and to %the infinite path $P_\infty$, which is
a connected graph $P_\infty$
having one vertex $v$ of degree 1 and $P_\infty -\{v\}$ with countable vertex set
and its vertices of degree 2 in $P_\infty$. If $r=1$, ($r=\infty$), we say that $P_r$
is a path of length 0, ($\infty$).

\begin{prop}
Each connected component of the subgraph ${\Lambda}[S]$ induced
by a {\rm PDS} $S$ in ${\Lambda}$ is a product subgraph
$P_r\times P_s$.
\end{prop}

\proof Assume that three vertices $u,v,w$ of
$Q=\{(i,j),(i+1,j),(i,j+1),(i+1,j+1)\}$ are members of a PDS $S$
in ${\Lambda}$. Clearly, $u,v,w$ are in a common component $C$
of ${\Lambda}[S]$. Then $Q\setminus\{u,v,w\}\subset C$
because of the definition of a PDS. By repeating
this argument, the conclusion follows. \qfd

Let ${\Gamma}_{m,n}=P_m\times P_n\subseteq{\Lambda}$ be the Cartesian product of
$P_m$ and $P_n$,
where $1\le m <\infty$ and $1\le n\le\infty$.
A horizontal, (vertical), {\it side} of ${\Gamma}_{m,n}$ is a
subgraph $P_m\times \{u\}$, ($\{u\}\times P_n$), where $u$ is an
endvertex of $P_n$, ($P_m$). If $1\le n\in\ZZ$, then
${\Gamma}_{m,n}$ is said to be an $m\times n$ {\it grid graph}. If
$n=1$, then ${\Gamma}_{m,n}$ is a path of length $m$. If
$n=\infty$, then ${\Gamma}_{m,n}={\Gamma}_{m,\infty}$,
an extended grid graph, is said to be a {\it band graph}, having just
one side of length $m-1$ and two infinite parallel sides.
The {\it width}, ({\it height}), of
${\Gamma}_{m,n}$ is the length of $P_m$, ($P_n$), that is $m-1$, ($n-1$).

\begin{cor}
Each component of the subgraph ${\Gamma}_{m,n}[S]$ induced by
a {\rm PDS} $S$ in ${\Gamma}_{m,n}$, where $0\leq m$ and
either $0\leq n\in\ZZ$, or $n=\infty$, is of the form $P_r\times
P_s$, where $1\leq r\leq m$ and either $1\leq s\leq n$ or
$s=\infty$. \qfd
\end{cor}

Without loss of generality, we identify
$V({\Gamma}_{m,n})$ with the set of vertices $(i,j)$ of ${\Lambda}$
with $0\leq i<m$, $0\leq j$ and if $n$ is finite then $j<n$. For $0\leq
j<n$, let the $j$-{\it level}
$H_j$ of ${\Gamma}_{m,n}$  be composed by those vertices of
${\Gamma}_{m,n}$ whose second coordinate is $j$. From now on,
Question 1.3 is considered with: {\bf(1)} $G_n={\Gamma}_{m,n}$,
for fixed $m>1$, where $0<n\in\ZZ$;
{\bf(2)} $G={\Gamma}_{m,\infty}$;
{\bf(3)} $G'=$
$H_0={\Gamma}_{m,1}=P_m\times\{0\}$.
In this formulation, we are asking for the existence of a PDS
$S=S(m,n,S')$ in ${\Gamma}_{m,n}$, for some $n>1$, with $S\cap G'=S'$.

\begin{prop}
For any $1<m\in\ZZ$, the band graph ${\Gamma}_{m,\infty}$ is
a periodic graph.
\end{prop}

\proof Set $G_n={\Gamma}_{m,n}$, for each $0<n\in\ZZ$, and set each $j$-level
$H_j=G_{j+1}-V(G_j)$. Then the union
${\Gamma}_{m,\infty}$ of the graphs ${\Gamma}_{m,n}$ is
periodic, as in the paragraph previous to Theorem 1.4. \qfd

If $S'\in\{\emptyset,H_0\}$ or otherwise if $S'$ dominates $H_0\setminus S'$, then
Question 1.3 is easily answered. In both cases, $S$ is
the Cartesian product of a path P times a PDS contained in another path $P'$,
or viceversa,
where ${\Gamma}_{m,n}$ is either $P\times P'$ or $P'\times P$, respectively.
Thus,
an admissible vertex subset $S'$ of $ H_0$ is said to be {\it
incomplete} if it is neither all of $H_0$ nor has all the vertices
of $H_0\setminus S'$ as dominated vertices, and
in considering Question 1.3, we require  from now on that
$S'$ be an incomplete admissible vertex subset, or IAVS, of $H_0$.

\begin{thm}
Let $0\le m\in\ZZ$ and let $S'$ be an IAVS of $G'$.
Then, there exists a continuation algorithm $\Theta$ that produces all sets $S(m,n,S')$, for
$0<n$ or $n=\infty$ and depends on eventual binary decisions,
with one of the two binary options being specifically greedy.
\end{thm}

\subsection{Plan of the remaining sections}

The algorithm $\Theta$ of Theorem 1.8 is discussed in Section 2 and
its running time, $O(2^{m+n})$, is attained in Section 3, while
Section 4 resets Question 1.3 under the scope of $\Theta$ and
considers its greediest strategy. In the extension of $\Theta$ to
band graphs of Subsection 4.2, Theorem 1.4 is invoked in order to
prune the binary decision tree of $\Theta$ into a finite threaded
tree, in which a closed walk is found that yields all possible
PDS-slices starting at $S'$, while in Section 5, the graphs induced
by the complements of PDSs in $G$ are codified by arrays of ordered
positive-integer pairs, via $\Theta$. In Section 6, a modification
of $\Theta$ attains all the grid graphs having total perfect codes,
i.e. PDSs with just 1-cubes as induced components, characterized by
Klostermeyer and Goldwasser, \cite{KG}. which allows to show that
there is only one total perfect code in ${\Lambda}$ that extends
them (and whose complement yields an aperiodic tiling, like the
Penrose tiling of \cite{Penrose}), despite their uncountability in
${\Lambda}$, found in Theorem 1 of \cite{DD1}.

\section{Continuation algorithm}

\newpage

We denote the algorithm $\Theta$
of Theorem 1.8 by $\Theta(\alpha,\beta)$, where
$\alpha$, ($\beta$), is a greedy, (non-greedy), option in the running of $\Theta$,
(see Subsection 2.3 below), so that
there is a binary decision tree $T_\Theta=T_{\Theta(\alpha,\beta)}$.
The root of $T_\Theta$ stands for the initial section of the running of $\Theta$ up to
the first binary decision instance.
Each internal node of $T_\Theta$ stands for a maximal section of the running of
$\Theta$ between two contiguous binary decision instances.

We agree that each node of $T_\Theta$, standing
for a partial solution to the problem of determining an
$S(m,n,S')$,  has a descending edge to
the left for option $\alpha$ and a descending edge to the right for option
$\beta$.

A {\it strategy}
$\sigma$ of $\Theta$ is an $\{\alpha,\beta\}$-sequence of
decisions determining a path from the root of $T_\Theta$, either
infinite, leading to a PDS $S(m,\infty,S')$ in a band graph
${\Gamma}_{m,\infty}$, or ending
at a leaf representing a PDS $S(m,n,S')$ in some grid graph ${\Gamma}_{m,n}$,
for $n<\infty$.

Algorithm $\Theta(\alpha,\beta)$, presented in the following four subsections,
is conceived with the second coordinate $j$
of vertices $(i,j)$ of ${\Gamma}_{m,n}$
increasing one unit per step, therefore advancing one level at a time, from $H_j$ to $H_{j+1}$,
in which candidate vertices for PDSs in ${\Gamma}_{m,n}$ or
${\Gamma}_{m,\infty}$ are considered.%

\subsection{Vertex-labeling initialization}

A function $f$ from $V({\Gamma}_{m,n})$
onto the alphabet $\{0,1,2,3,4\}$ is declared. We initialize $f$ by $f(v):=0$,
for every vertex $v$ of ${\Gamma}_{m,n}$ unless $j=0$ and either
$v\in V(H_0)\setminus S'$
is adjacent on its right, (left), to the leftmost, (rightmost),
vertex of a component of $ H_0[S']$, in which case we set $f(v):=1$, ($f(v):=3$), or
$v\in S'$, in which case we set $f(v):=2$.

\begin{prop}
An admissible vertex subset $S'$ of $ H_0$ is incomplete if and
only if $f(v)=0$ for some vertex $v$ of $ H_0$.
\end{prop}

\proof $f(v)=0$ for some vertex $v$ of $ H_0$ if and only if
either $v$ is an endvertex of $ H_0$ at distance $>1$ from $S'$ or
there exist two contiguous components of $ H_0[S']$ at distance
$>3$ such that a path between them contains $v$. \qfd

The algorithm $\Theta$ has an iterative body formed by a sequence
of steps, (presented in Subsections 2.2-4 below), that is applied
initially for $j:=0$.
This will produce a redefinition of $f$ on $H_1$ using labels $f(v)\in[0,4].$
Assuming that this
was already done on every $ H_\ell$, for
$0<\ell\leq j<n$, a further application of the iterative body of $\Theta$
redefines $f$ on $ H_{j+1}$.
In this case, we anticipate that: {\bf(a)} for any
vertex $w=(i,j+1)$ such that $\Theta$ cannot establish $f(w)\in[0,3]$, it is seen that
$f(i,j)=2$,
in which case $\Theta$ sets $f(w)=f(i,j+1)=4$; {\bf(b)}
if a vertex $u$ of $ H_{j+1}$ remains with $f(u)=0$,
then $\Theta$ cannot produce a PDS in
${\Gamma}_{m,j+2}$ and a further application of the body of $\Theta$ is needed;
in this case, the set
$W$ of vertices $w$ of ${\Gamma}_{m,j+2}$ with $f(w)=2$ will be
called a {\it quasiperfect dominating set}, or QPDS, in ${\mathcal
G}_{m,j+2}$, and $u$ is not dominated by $W$.

The five steps
contained in the following three subsections constitute the iterative body of $\Theta$,
until it stops in a passage of step 5 of Subsection 2.4.
Assume this body has run already on $H_\ell$, for $0\leq\ell\leq j$, so that $f$ is already redefined
on those levels.

\subsection{Labeling progressively the vertices of $ H_{j+1}$}

\begin{enumerate}
\item[1.] For $i:=0$ to $m-1$ do: If $f(i,j)=0$ then
\begin{enumerate}
\item[(A)] If $i>0$ and $f(i-1,j)>0$ then:
\begin{enumerate}\item[(a)] $f(i-1,j+1):=1$; \item[(b)]
If $i>1$ then:
\begin{enumerate}\item[(b1)] $k:=i-2$; \item[(b2)] While $f(k,j)=2$ do:

{(i)} $f(k,j+1):=4$;

{(ii)} $k:=k-1$;\end{enumerate}
\end{enumerate}
\item[(B)] If $i<m-1$ and $f(i+1,j)>0$ then:
\begin{enumerate}\item[(a)] $f(i+1,j+1):=3$;
\item[(b)] If $i<m-2$ then:
\begin{enumerate}\item[(b1)] $k:=i+2$; \item[(b2)]
While $f(k,j)=2$ do:

{(i)} $f(k,j+1):=4$;

{(ii)} $k:=k+1$;
\end{enumerate}\end{enumerate}
\item[(C)] $f(i,j+1):=2$.
\end{enumerate}
\item[2.] For $i:=0$ to $m-3$ do:

If $f(i+k,j)=k+1$ and $f(i+k,j+1)=0$, for $k=0,1,2$,
then:

$\,\,\,\,\,\,\,\,\,\,\,\,$ For $k:=0$ to 2 do: $f(i+k,j+1):=k+1$.
\end{enumerate}

\subsection{Binary decision instances}

\begin{enumerate}
\item[3.] If $f(0,j)=2$ and $f(0,j+1)=0$ then:
\begin{enumerate}\item[(a)] $k:=0$; \item[(b)]
While $f(k,j)=2$ and $f(k,j+1)=0$ do: $k:=k+1$; \item[(c)]
Select either option ($\alpha$) or option ($\beta$), where:
\begin{enumerate}
\item[($\alpha$)]

{(i)} For $\ell:=0$ to $k-1$ do: $f(\ell,j+1):=2$;

{(ii)} $f(k,j+1):=3$;
\item[($\beta$)] For $\ell:=0$ to $k-1$ do: $f(\ell,j+1):=4$;
\end{enumerate}\end{enumerate}\end{enumerate}
\begin{enumerate}         %1
\item[4.] For $i:=0$ to $m-2$ do: If $f(i,j)=1$ and $f(i+1,j+1)=0$
then:
\begin{enumerate}         %2
\item[(a)] $k:=i+1$; \item[(b)] While $f(k,j)=2$ and $f(k,j+1)=0$ do:
$k:=k+1$;
\item[(c)] If $k\leq m-1$ then:
If $f(k,j+1)=0$ then select either option ($\alpha$) or option ($\beta$), where:
\begin{enumerate}\item[]  %3
\begin{enumerate}         %4
\item[($\alpha$)]
(i) $f(i,j+1):=1$;

(ii) $f(k,j+1):=3$;

(iii) For $\ell:=i+1$ to $k-1$ do: $f(\ell,j+1):=2$;

\item[($\beta$)] For $\ell:=i+1$ to $k-1$ do: $f(\ell,j+1):=4$;

\end{enumerate}          %2
\end{enumerate}          %2
\noindent Else, (i.e. if $k=m$),
select either option ($\alpha$) or option ($\beta$), where:
\begin{enumerate}       %3
\item[]\begin{enumerate}%4
\item[($\alpha$)]
(i) $f(i,j+1):=1$;

(ii) For $\ell:=i+1$ to $m-1$ do: $f(\ell,j+1):=2$;
\item[($\beta$)] For $\ell:=i+1$ to $m-1$ do: $f(\ell,j+1):=4$.
\end{enumerate}         %3
\end{enumerate}         %2
\end{enumerate}         %1
\end{enumerate}         %0

\noindent
(Each decision taken by $\Theta(\alpha,\beta)$ involves either
option $\alpha$ setting some consecutive values $f(\ell,j+1)$ in the subset $\{1,2,3\}$, or
option $\beta$ setting those same values in the subset $\{0,4\}$, in such a way that $\beta$
leaves the leftmost and rightmost
values both as 0, while $\alpha$ leaves them greedily as 1 and 3, respectively,
so we may say that $\alpha$
is a {\it greedy} option and $\beta$ is not.
The decisions between options $\alpha$ and $\beta$ in
item 3(c) and the end of item 4(c) will be
referred to as binary {\it outer} decisions, or BOD. The initial
decision in item 4(c) will be referred to as a binary {\it inner}
decision or, BID).

\subsection{Checking PDS formation}

\newpage

\begin{enumerate}
\item[5.] Let $\tau(j)=|\{i\in[0,m):f(i,j+1)=0\}|$;

If $\tau(j)=0$, then the vertices $v=(i,\ell)$ for which
$0\leq\ell\leq j+1 \mbox{ and } f(i,\ell)=2$ constitute a PDS
$S$ in ${\Gamma}_{m,j+2}$;

Else: If $j<n$, then:
\begin{enumerate}
\item[{(a)}] $j:=j+1$; \item[{(b)}] go to Step
1.\end{enumerate}
\end{enumerate}

\begin{figure}
\unitlength=0.50mm
\special{em:linewidth 0.4pt}
\linethickness{0.4pt}
\begin{picture}(259.32,122.22)
\put(38.02,105.00){\circle{2.00}}
\put(48.02,105.00){\circle*{2.00}}
\put(58.02,105.00){\circle*{2.00}}
\put(68.02,105.00){\circle*{2.00}}
\put(78.02,105.00){\circle{2.00}}
\put(88.02,105.00){\circle{2.00}}
\put(98.02,105.00){\circle{2.00}}
\put(108.02,105.00){\circle{2.00}}
\put(118.02,105.00){\circle{2.00}}
\put(128.02,105.00){\circle*{2.00}}
\put(138.02,105.00){\circle{2.00}}
\put(148.02,105.00){\circle{2.00}}
\put(158.02,105.00){\circle{2.00}}
\put(168.02,105.00){\circle*{2.00}}
\put(178.02,105.00){\circle*{2.00}}
\put(188.02,105.00){\circle{2.00}}
\
\put(38.02,95.00){\circle{2.00}}
\put(48.02,95.00){\circle{2.00}}
\put(58.02,95.00){\circle{2.00}}
\put(68.02,95.00){\circle{2.00}}
\put(78.02,95.00){\circle{2.00}}
\put(88.02,95.00){\circle*{2.00}}
\put(98.02,95.00){\circle*{2.00}}
\put(108.02,95.00){\circle*{2.00}}
\put(118.02,95.00){\circle{2.00}}
\put(128.02,95.00){\circle{2.00}}
\put(138.02,95.00){\circle{2.00}}
\put(148.02,95.00){\circle*{2.00}}
\put(158.02,95.00){\circle{2.00}}
\put(168.02,95.00){\circle{2.00}}
\put(178.02,95.00){\circle{2.00}}
\put(188.02,95.00){\circle{2.00}}
\
\put(38.02,85.00){\circle*{2.00}}
\put(48.02,85.00){\circle{2.00}}
\put(58.02,85.00){\circle{2.00}}
\put(68.02,85.00){\circle{2.00}}
\put(78.02,85.00){\circle{2.00}}
\put(88.02,85.00){\circle{2.00}}
\put(98.02,85.00){\circle{2.00}}
\put(108.02,85.00){\circle{2.00}}
\put(118.02,85.00){\circle{2.00}}
\put(128.02,85.00){\circle{2.00}}
\put(138.02,85.00){\circle{2.00}}
\put(148.02,85.00){\circle*{2.00}}
\put(158.02,85.00){\circle{2.00}}
\put(168.02,85.00){\circle{2.00}}
\put(178.02,85.00){\circle{2.00}}
\put(188.02,85.00){\circle*{2.00}}
\
\put(38.02,75.00){\circle{2.00}}
\put(48.02,75.00){\circle{2.00}}
\put(58.02,75.00){\circle*{2.00}}
\put(68.02,75.00){\circle*{2.00}}
\put(78.02,75.00){\circle*{2.00}}
\put(88.02,75.00){\circle{2.00}}
\put(98.02,75.00){\circle{2.00}}
\put(108.02,75.00){\circle{2.00}}
\put(118.02,75.00){\circle*{2.00}}
\put(128.02,75.00){\circle*{2.00}}
\put(138.02,75.00){\circle{2.00}}
\put(148.02,75.00){\circle{2.00}}
\put(158.02,75.00){\circle{2.00}}
\put(168.02,75.00){\circle*{2.00}}
\put(178.02,75.00){\circle{2.00}}
\put(188.02,75.00){\circle{2.00}}
\
\put(38.02,65.00){\circle{2.00}}
\put(78.02,65.00){\circle{2.00}}
\put(88.02,65.00){\circle{2.00}}
\put(98.02,65.00){\circle*{2.00}}
\put(108.02,65.00){\circle{2.00}}
\put(118.02,65.00){\circle{2.00}}
\put(128.02,65.00){\circle{2.00}}
\put(138.02,65.00){\circle{2.00}}
\put(148.02,65.00){\circle{2.00}}
\put(158.02,65.00){\circle{2.00}}
\put(168.02,65.00){\circle*{2.00}}
\put(178.02,65.00){\circle{2.00}}
\put(188.02,65.00){\circle{2.00}}
\
\put(38.02,55.00){\circle*{2.00}}
\put(48.02,55.00){\circle*{2.00}}
\put(58.02,55.00){\circle{2.00}}
\put(68.02,55.00){\circle{2.00}}
\put(78.02,55.00){\circle{2.00}}
\put(88.02,55.00){\circle{2.00}}
\put(98.02,55.00){\circle*{2.00}}
\put(108.02,55.00){\circle{2.00}}
\put(118.02,55.00){\circle{2.00}}
\put(128.02,55.00){\circle{2.00}}
\put(138.02,55.00){\circle*{2.00}}
\put(148.02,55.00){\circle*{2.00}}
\put(158.02,55.00){\circle{2.00}}
\put(168.02,55.00){\circle{2.00}}
\put(178.02,55.00){\circle{2.00}}
\put(188.02,55.00){\circle*{2.00}}
\
\put(38.02,45.00){\circle{2.00}}
\put(48.02,45.00){\circle{2.00}}
\put(58.02,45.00){\circle{2.00}}
\put(68.02,45.00){\circle*{2.00}}
\put(78.02,45.00){\circle*{2.00}}
\put(88.02,45.00){\circle{2.00}}
\put(98.02,45.00){\circle{2.00}}
\put(108.02,45.00){\circle{2.00}}
\put(118.02,45.00){\circle*{2.00}}
\put(128.02,45.00){\circle{2.00}}
\put(138.02,45.00){\circle{2.00}}
\put(148.02,45.00){\circle{2.00}}
\put(158.02,45.00){\circle{2.00}}
\put(168.02,45.00){\circle{2.00}}
\put(178.02,45.00){\circle{2.00}}
\put(188.02,45.00){\circle*{2.00}}
\
\put(38.02,35.00){\circle{2.00}}
\put(48.02,35.00){\circle{2.00}}
\put(58.02,35.00){\circle{2.00}}
\put(68.02,35.00){\circle{2.00}}
\put(78.02,35.00){\circle{2.00}}
\put(88.02,35.00){\circle{2.00}}
\put(98.02,35.00){\circle{2.00}}
\put(108.02,35.00){\circle{2.00}}
\put(118.02,35.00){\circle*{2.00}}
\put(128.02,35.00){\circle{2.00}}
\put(138.02,35.00){\circle{2.00}}
\put(148.02,35.00){\circle{2.00}}
\
\put(38.02,25.00){\circle*{2.00}}
\put(48.02,25.00){\circle*{2.00}}
\put(58.02,25.00){\circle*{2.00}}
\put(68.02,25.00){\circle{2.00}}
\put(78.02,25.00){\circle{2.00}}
\put(88.02,25.00){\circle*{2.00}}
\put(98.02,25.00){\circle*{2.00}}
\put(108.02,25.00){\circle{2.00}}
\put(118.02,25.00){\circle{2.00}}
\put(128.02,25.00){\circle{2.00}}
\put(138.02,25.00){\circle*{2.00}}
\put(148.02,25.00){\circle{2.00}}
\put(158.02,25.00){\circle{2.00}}
\put(168.02,25.00){\circle{2.00}}
\put(178.02,25.00){\circle{2.00}}
\put(188.02,25.00){\circle{2.00}}
\
\put(38.02,15.00){\circle*{2.00}}
\put(48.02,15.00){\circle*{2.00}}
\put(58.02,15.00){\circle*{2.00}}
\put(68.02,15.00){\circle{2.00}}
\put(78.02,15.00){\circle{2.00}}
\put(88.02,15.00){\circle{2.00}}
\put(98.02,15.00){\circle{2.00}}
\put(108.02,15.00){\circle{2.00}}
\put(118.02,15.00){\circle{2.00}}
\put(128.02,15.00){\circle{2.00}}
\put(138.02,15.00){\circle*{2.00}}
\put(148.02,15.00){\circle{2.00}}
\put(158.02,15.00){\circle{2.00}}
\put(168.02,15.00){\circle{2.00}}
\put(178.02,15.00){\circle*{2.00}}
\put(188.02,15.00){\circle*{2.00}}
\
\put(38.02,5.00){\circle{2.00}}
\put(48.02,5.00){\circle{2.00}}
\put(58.02,5.00){\circle{2.00}}
\put(68.02,5.00){\circle{2.00}}
\put(78.02,5.00){\circle*{2.00}}
\put(88.02,5.00){\circle{2.00}}
\put(98.02,5.00){\circle{2.00}}
\put(108.02,5.00){\circle*{2.00}}
\put(118.02,5.00){\circle*{2.00}}
\put(128.02,5.00){\circle{2.00}}
\put(138.02,5.00){\circle{2.00}}
\put(148.02,5.00){\circle{2.00}}
\put(158.02,5.00){\circle*{2.00}}
\put(168.02,5.00){\circle{2.00}}
\put(178.02,5.00){\circle{2.00}}
\put(188.02,5.00){\circle{2.00}}
\put(48.02,65.00){\circle{2.00}}
\put(58.02,65.00){\circle{2.00}}
\put(68.02,65.00){\circle{2.00}}
\put(158.02,35.00){\circle*{2.00}}
\put(168.02,35.00){\circle*{2.00}}
\put(178.02,35.00){\circle{2.00}}
\put(188.02,35.00){\circle{2.00}}
\put(38.02,108.00){\makebox(0,0)[cc]{$_1$}}
\put(48.02,108.00){\makebox(0,0)[cc]{$_2$}}
\put(58.02,108.00){\makebox(0,0)[cc]{$_2$}}
\put(68.02,108.00){\makebox(0,0)[cc]{$_2$}}
\put(78.02,108.00){\makebox(0,0)[cc]{$_3$}}
\put(88.02,108.00){\makebox(0,0)[cc]{$_0$}}
\put(98.02,108.00){\makebox(0,0)[cc]{$_0$}}
\put(108.02,108.00){\makebox(0,0)[cc]{$_0$}}
\put(118.02,108.00){\makebox(0,0)[cc]{$_1$}}
\put(128.02,108.00){\makebox(0,0)[cc]{$_2$}}
\put(138.02,108.00){\makebox(0,0)[cc]{$_3$}}
\put(148.02,108.00){\makebox(0,0)[cc]{$_0$}}
\put(158.02,108.00){\makebox(0,0)[cc]{$_1$}}
\put(168.02,108.00){\makebox(0,0)[cc]{$_2$}}
\put(178.02,108.00){\makebox(0,0)[cc]{$_2$}}
\put(188.02,108.00){\makebox(0,0)[cc]{$_3$}}
\put(38.02,98.00){\makebox(0,0)[cc]{$_0$}}
\put(48.02,98.00){\makebox(0,0)[cc]{$_4$}}
\put(58.02,98.00){\makebox(0,0)[cc]{$_4$}}
\put(68.02,98.00){\makebox(0,0)[cc]{$_4$}}
\put(78.02,98.00){\makebox(0,0)[cc]{$_1$}}
\put(88.02,98.00){\makebox(0,0)[cc]{$_2$}}
\put(98.02,98.00){\makebox(0,0)[cc]{$_2$}}
\put(108.02,98.00){\makebox(0,0)[cc]{$_2$}}
\put(118.02,98.00){\makebox(0,0)[cc]{$_3$}}
\put(128.02,98.00){\makebox(0,0)[cc]{$_4$}}
\put(138.02,98.00){\makebox(0,0)[cc]{$_1$}}
\put(148.02,98.00){\makebox(0,0)[cc]{$_2$}}
\put(158.02,98.00){\makebox(0,0)[cc]{$_3$}}
\put(168.02,98.00){\makebox(0,0)[cc]{$_4$}}
\put(178.02,98.00){\makebox(0,0)[cc]{$_4$}}
\put(188.02,98.00){\makebox(0,0)[cc]{$_0$}}
\put(38.02,88.00){\makebox(0,0)[cc]{$_2$}}
\put(48.02,88.00){\makebox(0,0)[cc]{$_3$}}
\put(58.02,88.00){\makebox(0,0)[cc]{$_0$}}
\put(68.02,88.00){\makebox(0,0)[cc]{$_0$}}
\put(78.02,88.00){\makebox(0,0)[cc]{$_0$}}
\put(88.02,88.00){\makebox(0,0)[cc]{$_4$}}
\put(98.02,88.00){\makebox(0,0)[cc]{$_4$}}
\put(108.02,88.00){\makebox(0,0)[cc]{$_4$}}
\put(118.02,88.00){\makebox(0,0)[cc]{$_0$}}
\put(128.02,88.00){\makebox(0,0)[cc]{$_0$}}
\put(138.02,88.00){\makebox(0,0)[cc]{$_1$}}
\put(148.02,88.00){\makebox(0,0)[cc]{$_2$}}
\put(158.02,88.00){\makebox(0,0)[cc]{$_3$}}
\put(168.02,88.00){\makebox(0,0)[cc]{$_0$}}
\put(178.02,88.00){\makebox(0,0)[cc]{$_1$}}
\put(188.02,88.00){\makebox(0,0)[cc]{$_2$}}
\put(38.02,78.00){\makebox(0,0)[cc]{$_4$}}
\put(48.02,78.00){\makebox(0,0)[cc]{$_1$}}
\put(58.02,78.00){\makebox(0,0)[cc]{$_2$}}
\put(68.02,78.00){\makebox(0,0)[cc]{$_2$}}
\put(78.02,78.00){\makebox(0,0)[cc]{$_2$}}
\put(88.02,78.00){\makebox(0,0)[cc]{$_3$}}
\put(98.02,78.00){\makebox(0,0)[cc]{$_0$}}
\put(108.02,78.00){\makebox(0,0)[cc]{$_1$}}
\put(118.02,78.00){\makebox(0,0)[cc]{$_2$}}
\put(128.02,78.00){\makebox(0,0)[cc]{$_2$}}
\put(138.02,78.00){\makebox(0,0)[cc]{$_3$}}
\put(148.02,78.00){\makebox(0,0)[cc]{$_4$}}
\put(158.02,78.00){\makebox(0,0)[cc]{$_1$}}
\put(168.02,78.00){\makebox(0,0)[cc]{$_2$}}
\put(178.02,78.00){\makebox(0,0)[cc]{$_3$}}
\put(188.02,78.00){\makebox(0,0)[cc]{$_4$}}
\put(38.02,68.00){\makebox(0,0)[cc]{$_0$}}
\put(48.02,68.00){\makebox(0,0)[cc]{$_0$}}
\put(58.02,68.00){\makebox(0,0)[cc]{$_4$}}
\put(68.02,68.00){\makebox(0,0)[cc]{$_4$}}
\put(78.02,68.00){\makebox(0,0)[cc]{$_4$}}
\put(88.02,68.00){\makebox(0,0)[cc]{$_1$}}
\put(98.02,68.00){\makebox(0,0)[cc]{$_2$}}
\put(108.02,68.00){\makebox(0,0)[cc]{$_3$}}
\put(118.02,68.00){\makebox(0,0)[cc]{$_4$}}
\put(128.02,68.00){\makebox(0,0)[cc]{$_4$}}
\put(138.02,68.00){\makebox(0,0)[cc]{$_0$}}
\put(148.02,68.00){\makebox(0,0)[cc]{$_0$}}
\put(156.02,68.00){\makebox(0,0)[cc]{$_1$}}
\put(168.02,68.00){\makebox(0,0)[cc]{$_2$}}
\put(176.02,68.00){\makebox(0,0)[cc]{$_3$}}
\put(188.02,68.00){\makebox(0,0)[cc]{$_0$}}
\put(38.02,58.00){\makebox(0,0)[cc]{$_2$}}
\put(48.02,58.00){\makebox(0,0)[cc]{$_2$}}
\put(58.02,58.00){\makebox(0,0)[cc]{$_3$}}
\put(68.02,58.00){\makebox(0,0)[cc]{$_0$}}
\put(78.02,58.00){\makebox(0,0)[cc]{$_0$}}
\put(86.02,58.00){\makebox(0,0)[cc]{$_1$}}
\put(98.02,58.00){\makebox(0,0)[cc]{$_2$}}
\put(106.02,58.00){\makebox(0,0)[cc]{$_3$}}
\put(118.02,58.00){\makebox(0,0)[cc]{$_0$}}
\put(128.02,58.00){\makebox(0,0)[cc]{$_1$}}
\put(138.02,58.00){\makebox(0,0)[cc]{$_2$}}
\put(148.02,58.00){\makebox(0,0)[cc]{$_2$}}
\put(158.02,58.00){\makebox(0,0)[cc]{$_3$}}
\put(168.02,58.00){\makebox(0,0)[cc]{$_4$}}
\put(178.02,58.00){\makebox(0,0)[cc]{$_1$}}
\put(188.02,58.00){\makebox(0,0)[cc]{$_2$}}
\put(38.02,48.00){\makebox(0,0)[cc]{$_4$}}
\put(48.02,48.00){\makebox(0,0)[cc]{$_4$}}
\put(58.02,48.00){\makebox(0,0)[cc]{$_1$}}
\put(68.02,48.00){\makebox(0,0)[cc]{$_2$}}
\put(78.02,48.00){\makebox(0,0)[cc]{$_2$}}
\put(88.02,48.00){\makebox(0,0)[cc]{$_3$}}
\put(98.02,48.00){\makebox(0,0)[cc]{$_4$}}
\put(108.02,48.00){\makebox(0,0)[cc]{$_1$}}
\put(118.02,48.00){\makebox(0,0)[cc]{$_2$}}
\put(128.02,48.00){\makebox(0,0)[cc]{$_3$}}
\put(138.02,48.00){\makebox(0,0)[cc]{$_4$}}
\put(148.02,48.00){\makebox(0,0)[cc]{$_4$}}
\put(158.02,48.00){\makebox(0,0)[cc]{$_0$}}
\put(168.02,48.00){\makebox(0,0)[cc]{$_0$}}
\put(178.02,48.00){\makebox(0,0)[cc]{$_1$}}
\put(188.02,48.00){\makebox(0,0)[cc]{$_2$}}
\put(38.02,38.00){\makebox(0,0)[cc]{$_0$}}
\put(48.02,38.00){\makebox(0,0)[cc]{$_0$}}
\put(58.02,38.00){\makebox(0,0)[cc]{$_0$}}
\put(68.02,38.00){\makebox(0,0)[cc]{$_4$}}
\put(78.02,38.00){\makebox(0,0)[cc]{$_4$}}
\put(88.02,38.00){\makebox(0,0)[cc]{$_0$}}
\put(98.02,38.00){\makebox(0,0)[cc]{$_0$}}
\put(106.02,38.00){\makebox(0,0)[cc]{$_1$}}
\put(118.02,38.00){\makebox(0,0)[cc]{$_2$}}
\put(128.02,38.00){\makebox(0,0)[cc]{$_3$}}
\put(138.02,38.00){\makebox(0,0)[cc]{$_0$}}
\put(148.02,38.00){\makebox(0,0)[cc]{$_1$}}
\put(158.02,38.00){\makebox(0,0)[cc]{$_2$}}
\put(168.02,38.00){\makebox(0,0)[cc]{$_2$}}
\put(178.02,38.00){\makebox(0,0)[cc]{$_3$}}
\put(188.02,38.00){\makebox(0,0)[cc]{$_4$}}
\put(38.02,28.00){\makebox(0,0)[cc]{$_2$}}
\put(48.02,28.00){\makebox(0,0)[cc]{$_2$}}
\put(58.02,28.00){\makebox(0,0)[cc]{$_2$}}
\put(68.02,28.00){\makebox(0,0)[cc]{$_3$}}
\put(78.02,28.00){\makebox(0,0)[cc]{$_1$}}
\put(88.02,28.00){\makebox(0,0)[cc]{$_2$}}
\put(98.02,28.00){\makebox(0,0)[cc]{$_2$}}
\put(108.02,28.00){\makebox(0,0)[cc]{$_3$}}
\put(118.02,28.00){\makebox(0,0)[cc]{$_4$}}
\put(128.02,28.00){\makebox(0,0)[cc]{$_1$}}
\put(138.02,28.00){\makebox(0,0)[cc]{$_2$}}
\put(148.02,28.00){\makebox(0,0)[cc]{$_3$}}
\put(158.02,28.00){\makebox(0,0)[cc]{$_4$}}
\put(168.02,28.00){\makebox(0,0)[cc]{$_4$}}
\put(178.02,28.00){\makebox(0,0)[cc]{$_0$}}
\put(188.02,28.00){\makebox(0,0)[cc]{$_0$}}
\put(38.02,18.00){\makebox(0,0)[cc]{$_2$}}
\put(48.02,18.00){\makebox(0,0)[cc]{$_2$}}
\put(58.02,18.00){\makebox(0,0)[cc]{$_2$}}
\put(68.02,18.00){\makebox(0,0)[cc]{$_3$}}
\put(78.02,18.00){\makebox(0,0)[cc]{$_0$}}
\put(88.02,18.00){\makebox(0,0)[cc]{$_4$}}
\put(98.02,18.00){\makebox(0,0)[cc]{$_4$}}
\put(108.02,18.00){\makebox(0,0)[cc]{$_0$}}
\put(118.02,18.00){\makebox(0,0)[cc]{$_0$}}
\put(128.02,18.00){\makebox(0,0)[cc]{$_1$}}
\put(138.02,18.00){\makebox(0,0)[cc]{$_2$}}
\put(146.02,18.00){\makebox(0,0)[cc]{$_3$}}
\put(158.02,18.00){\makebox(0,0)[cc]{$_0$}}
\put(168.02,18.00){\makebox(0,0)[cc]{$_1$}}
\put(178.02,18.00){\makebox(0,0)[cc]{$_2$}}
\put(188.02,18.00){\makebox(0,0)[cc]{$_2$}}
\put(38.02,8.00){\makebox(0,0)[cc]{$_4$}}
\put(48.02,8.00){\makebox(0,0)[cc]{$_4$}}
\put(58.02,8.00){\makebox(0,0)[cc]{$_4$}}
\put(68.02,8.00){\makebox(0,0)[cc]{$_1$}}
\put(78.02,8.00){\makebox(0,0)[cc]{$_2$}}
\put(88.02,8.00){\makebox(0,0)[cc]{$_3$}}
\put(98.02,8.00){\makebox(0,0)[cc]{$_1$}}
\put(108.02,8.00){\makebox(0,0)[cc]{$_2$}}
\put(118.02,8.00){\makebox(0,0)[cc]{$_2$}}
\put(128.02,8.00){\makebox(0,0)[cc]{$_3$}}
\put(138.02,8.00){\makebox(0,0)[cc]{$_4$}}
\put(148.02,8.00){\makebox(0,0)[cc]{$_1$}}
\put(158.02,8.00){\makebox(0,0)[cc]{$_2$}}
\put(168.02,8.00){\makebox(0,0)[cc]{$_3$}}
\put(178.02,8.00){\makebox(0,0)[cc]{$_4$}}
\put(188.02,8.00){\makebox(0,0)[cc]{$_4$}}
\put(28.02,105.00){\makebox(0,0)[cc]{$^0$}}
\put(28.02,95.00){\makebox(0,0)[cc]{$^1$}}
\put(28.02,85.00){\makebox(0,0)[cc]{$^2$}}
\put(28.02,75.00){\makebox(0,0)[cc]{$^3$}}
\put(28.02,65.00){\makebox(0,0)[cc]{$^4$}}
\put(28.02,55.00){\makebox(0,0)[cc]{$^5$}}
\put(28.02,45.00){\makebox(0,0)[cc]{$^6$}}
\put(28.02,35.00){\makebox(0,0)[cc]{$^7$}}
\put(28.02,25.00){\makebox(0,0)[cc]{$^8$}}
\put(28.02,15.00){\makebox(0,0)[cc]{$^9$}}
\put(28.02,5.00){\makebox(0,0)[cc]{$^{10}$}}
\put(38.85,118.89){\makebox(0,0)[cc]{$^0$}}
\put(48.02,118.89){\makebox(0,0)[cc]{$^1$}}
\put(58.02,118.89){\makebox(0,0)[cc]{$^2$}}
\put(68.02,118.89){\makebox(0,0)[cc]{$^3$}}
\put(78.02,118.89){\makebox(0,0)[cc]{$^4$}}
\put(88.02,118.89){\makebox(0,0)[cc]{$^5$}}
\put(98.02,118.89){\makebox(0,0)[cc]{$^6$}}
\put(108.02,118.89){\makebox(0,0)[cc]{$^7$}}
\put(118.02,118.89){\makebox(0,0)[cc]{$^8$}}
\put(128.02,118.89){\makebox(0,0)[cc]{$^9$}}
\put(138.02,118.89){\makebox(0,0)[cc]{$^{10}$}}
\put(148.02,118.89){\makebox(0,0)[cc]{$^{11}$}}
\put(158.02,118.89){\makebox(0,0)[cc]{$^{12}$}}
\put(168.02,118.89){\makebox(0,0)[cc]{$^{13}$}}
\put(178.02,118.89){\makebox(0,0)[cc]{$^{14}$}}
\put(188.02,118.89){\makebox(0,0)[cc]{$^{15}$}}
\put(26.91,115.00){\makebox(0,0)[cc]{$_{j=}$}}
\put(33.02,120.00){\makebox(0,0)[cc]{$_{i=}$}}
\emline{30.91}{116.00}{1}{30.91}{116.00}{2}
\emline{30.91}{116.00}{3}{30.91}{0.00}{4}
\emline{31.08}{115.83}{5}{194.83}{115.83}{6}
\put(184.06,47.08){\oval(22.08,10.00)[]}
\put(98.85,87.08){\oval(50.00,10.00)[]}
\put(72.78,36.43){\oval(40.00,9.05)[]}
\put(52.78,16.43){\oval(40.00,9.05)[]}
\put(92.78,16.43){\oval(40.00,9.05)[]}
\emline{31.43}{116.00}{7}{25.21}{122.22}{8}
\put(216.02,105.00){\circle{2.00}}
\put(226.02,105.00){\circle*{2.00}}
\put(216.02,95.00){\circle{2.00}}
\put(226.02,95.00){\circle{2.00}}
\put(236.02,95.00){\circle{2.00}}
\put(216.02,85.00){\circle*{2.00}}
\put(226.02,85.00){\circle{2.00}}
\put(236.02,85.00){\circle{2.00}}
\put(246.02,85.00){\circle{2.00}}
\put(256.02,85.00){\circle{2.00}}
\put(216.02,75.00){\circle{2.00}}
\put(226.02,75.00){\circle{2.00}}
\put(236.02,75.00){\circle*{2.00}}
\put(216.02,65.00){\circle{2.00}}
\put(216.02,55.00){\circle*{2.00}}
\put(226.02,55.00){\circle*{2.00}}
\put(236.02,55.00){\circle{2.00}}
\put(246.02,55.00){\circle{2.00}}
\put(256.02,55.00){\circle{2.00}}
\put(216.02,45.00){\circle{2.00}}
\put(226.02,45.00){\circle{2.00}}
\put(236.02,45.00){\circle{2.00}}
\put(246.02,45.00){\circle*{2.00}}
\put(226.02,65.00){\circle{2.00}}
\put(236.02,65.00){\circle{2.00}}
\put(246.02,65.00){\circle{2.00}}
\put(216.02,108.00){\makebox(0,0)[cc]{$_1$}}
\put(226.02,108.00){\makebox(0,0)[cc]{$_2$}}
\put(236.02,108.00){\makebox(0,0)[cc]{$_3$}}
\put(246.02,108.00){\makebox(0,0)[cc]{$_0$}}
\put(256.02,108.00){\makebox(0,0)[cc]{$_0$}}
\put(216.02,98.00){\makebox(0,0)[cc]{$_0$}}
\put(226.02,98.00){\makebox(0,0)[cc]{$_4$}}
\put(236.02,98.00){\makebox(0,0)[cc]{$_1$}}
\put(246.02,98.00){\makebox(0,0)[cc]{$_2$}}
\put(256.02,98.00){\makebox(0,0)[cc]{$_2$}}
\put(216.02,88.00){\makebox(0,0)[cc]{$_2$}}
\put(226.02,88.00){\makebox(0,0)[cc]{$_3$}}
\put(236.02,88.00){\makebox(0,0)[cc]{$_0$}}
\put(246.02,88.00){\makebox(0,0)[cc]{$_4$}}
\put(256.02,88.00){\makebox(0,0)[cc]{$_4$}}
\put(216.02,78.00){\makebox(0,0)[cc]{$_4$}}
\put(226.02,78.00){\makebox(0,0)[cc]{$_1$}}
\put(236.02,78.00){\makebox(0,0)[cc]{$_2$}}
\put(246.02,78.00){\makebox(0,0)[cc]{$_3$}}
\put(256.02,78.00){\makebox(0,0)[cc]{$_0$}}
\put(216.02,68.00){\makebox(0,0)[cc]{$_0$}}
\put(226.02,68.00){\makebox(0,0)[cc]{$_0$}}
\put(236.02,68.00){\makebox(0,0)[cc]{$_4$}}
\put(246.02,68.00){\makebox(0,0)[cc]{$_1$}}
\put(256.02,68.00){\makebox(0,0)[cc]{$_2$}}
\put(216.02,58.00){\makebox(0,0)[cc]{$_2$}}
\put(226.02,58.00){\makebox(0,0)[cc]{$_2$}}
\put(236.02,58.00){\makebox(0,0)[cc]{$_3$}}
\put(246.02,58.00){\makebox(0,0)[cc]{$_0$}}
\put(256.02,58.00){\makebox(0,0)[cc]{$_4$}}
\put(216.02,48.00){\makebox(0,0)[cc]{$_4$}}
\put(226.02,48.00){\makebox(0,0)[cc]{$_4$}}
\put(236.02,48.00){\makebox(0,0)[cc]{$_1$}}
\put(246.02,48.00){\makebox(0,0)[cc]{$_2$}}
\put(256.02,48.00){\makebox(0,0)[cc]{$_3$}}
\put(216.85,118.89){\makebox(0,0)[cc]{$^0$}}
\put(226.02,118.89){\makebox(0,0)[cc]{$^1$}}
\put(236.02,118.89){\makebox(0,0)[cc]{$^2$}}
\put(246.02,118.89){\makebox(0,0)[cc]{$^3$}}
\put(256.02,118.89){\makebox(0,0)[cc]{$^4$}}
\emline{208.91}{116.00}{9}{208.91}{116.00}{10}
\emline{208.91}{116.00}{11}{208.91}{0.00}{12}
\put(236.02,105.00){\circle{2.00}}
\put(246.02,105.00){\circle{2.00}}
\put(256.02,105.00){\circle{2.00}}
\put(246.02,95.00){\circle*{2.00}}
\put(256.02,95.00){\circle*{2.00}}
\put(246.02,75.00){\circle{2.00}}
\put(256.02,75.00){\circle{2.00}}
\put(256.02,65.00){\circle*{2.00}}
\put(256.02,55.00){\circle{2.00}}
\put(256.02,45.00){\circle{2.00}}
\put(245.82,87.00){\oval(27.00,10.00)[]}
\put(250.82,57.00){\oval(17.00,10.00)[]}
\put(232.32,92.00){\makebox(0,0)[cc]{$_\beta$}}
\put(242.32,62.00){\makebox(0,0)[cc]{$_\beta$}}
\put(206.02,105.00){\makebox(0,0)[cc]{$^0$}}
\put(206.02,95.00){\makebox(0,0)[cc]{$^1$}}
\put(206.02,85.00){\makebox(0,0)[cc]{$^2$}}
\put(206.02,75.00){\makebox(0,0)[cc]{$^3$}}
\put(206.02,65.00){\makebox(0,0)[cc]{$^4$}}
\put(206.02,55.00){\makebox(0,0)[cc]{$^5$}}
\put(206.02,45.00){\makebox(0,0)[cc]{$^6$}}
\put(204.91,115.00){\makebox(0,0)[cc]{$_{j=}$}}
\put(211.02,120.00){\makebox(0,0)[cc]{$_{i=}$}}
\emline{209.43}{116.00}{13}{203.21}{122.22}{14}
\emline{208.21}{116.11}{15}{259.32}{116.11}{16}
\put(173.21,52.00){\makebox(0,0)[cc]{$_\alpha$}}
\put(74.21,92.00){\makebox(0,0)[cc]{$_\beta$}}
\put(54.21,42.00){\makebox(0,0)[cc]{$_\beta$}}
\put(34.21,22.00){\makebox(0,0)[cc]{$_\alpha$}}
\put(74.21,22.00){\makebox(0,0)[cc]{$_\beta$}}
\end{picture}
\caption{Two applications of $\Theta$}
\end{figure}
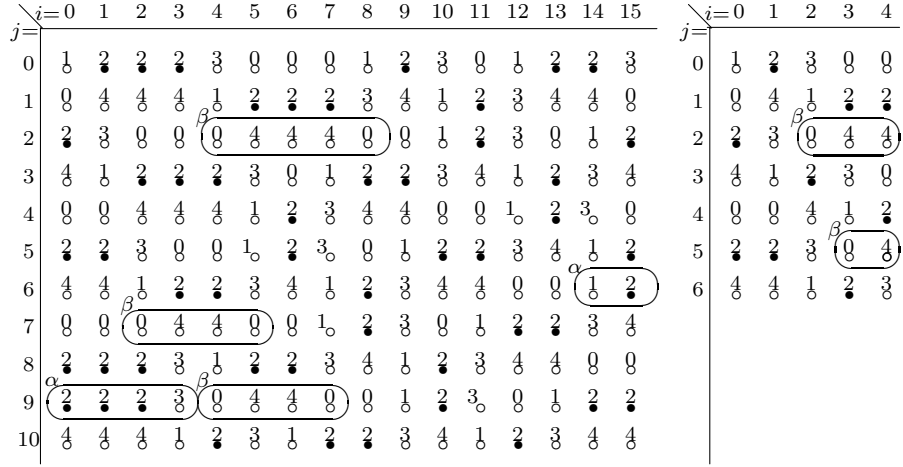

\subsection{Examples of the running of $\Theta(\alpha,\beta)$}

A resulting QPDS or PDS $S$ %in ${\Gamma}_{m,n}$
is given by means of black nodes. %and the remaining ones by white nodes.
Figure 1 depicts two PDSs
in ${\Gamma}_{m,n}$, both obtained via $\Theta$:
for $(m,n)=(16,11)$ on the left and for
$(m,n)=(5,7)$ on the right.
Those vertex subsets affected by binary decisions in
$\{\alpha,\beta\}$ are enclosed in ovals.
Observe we are not indicating the edges of ${\Gamma}_{m,n}$.

In the example of a PDS in ${\Gamma}_{16,11}$, ${\mathcal L}[S']$ is
formed by three components induced by $\{(1,0),(2,0),(3,0)\},
\{(9,0)\}$ and $\{(13,0),(14,0)\}.$ After initializing $f$, the
labeling of $H_0$ for $i\in[0,m)=[0,16)$ is
$f(H_0)=1222300012301223$. After Step 1, %(Subsection 2.2),
the labeling of $H_1$ is $f(H_1)$ $=0444122234123440$, where the 0's
are remnants of the initialization of $f$ in this level. Steps 2 to
5 %, (Subsections 2.3-4),
are passed through with no variations and
in the resulting first instance of Step 5, after setting $j=1$, the
process returns to Step 1 with $j=2$. After the new instance of Step
1, $f(H_2)$ looks like 2300?444?0???012, where `?' stands for those
values of $f$ equal to 0 that may still vary in the running of
$\Theta$. After the subsequent instance of Step 2, $f(H_2)$ looks
like 2300?444?0123012. A first decision must be taken now according
to step 4, for $i=4$ and $k=8$. Figure 1 shows that $\beta$ was
taken, yielding $f(H_2)=2300044400123012$. The rest of $\Theta$ is
applied similarly. Figure 1 shows that the successive decisions
happened at $j=$2,6,7,9,9 for respective steps
4(BID),4(BOD),4(BID),3(BOD),4(BID) with
$(i,k)=(4,8),k=14,(i,k)=(2,5),k=3,(i,k)=(4,6)$ and that the options
followed the strategy $\sigma=\{\beta,\alpha,\beta,\alpha,\beta\}$.
We leave for the reader to check the example at the right of Figure
1.

\section{PDS generation and exhaustion}

\begin{thm}
If, for some $0<j\in\ZZ$, the running of $\Theta$ for a
particular strategy $\sigma$ yields an instance of step $5$ of Subsection 2.4 with
$\tau(j)=0$, then it also yields, for the subsequence of $\sigma$
covered during the running, a {\rm PDS} $S(m,j+2,S')$ formed by those
vertices $v$ for which $f(v)=2$.
\end{thm}

\proof
We will say that a subsequence $\zeta$ of $f(H_j)$ is {\it
$\kappa$-covered}, where
$\kappa\in[0,4]$, if: {\bf (a)} the subgraph induced by the vertices
of $\zeta$ is a path; {\bf (b)} the value $f(i,j)$ associated to each
of the vertices of $\zeta$ is $\kappa$ and
{\bf (c)} $\zeta$ is maximal with respect to (a)-(b).
Assume the $f$-labeling at $ H_0, H_1,\ldots,H_j$ was
completed by $\Theta$, so that the running of $\Theta$ returns to step 1,
forcing new values of $f$ at level $ H_{j+1}$ with a widening
`bell' or `gable' effect: every component $ H_j[\{(i,j),\ldots,(k,j)\}]$ of
the subgraph induced by the 0-labelled vertices at level $ H_j$,
($0\leq i\leq k<m$), determines a corresponding component
$H_{j+1}[\{(i,j+1),\ldots,(k,j+1)\}]$ of the subgraph induced by
the new 2-labelled vertices at level $ H_{j+1}$, (step 1(C)), plus
a 1-labelled vertex at position $(i-1,j+1)$, if $k+1<m$, (step
1(A)(a)), and a 3-labelled-vertex at position $(k+1,j+1)$, if
$i\geq 1$, (step 1(B)(a)), as in the two upper levels
corresponding to each polygon of ${\Gamma}_{m,j+2}$ in Figure
1. Thus, a `bell island' or `pediment' given as a subsequence of contiguous
labels 122$\ldots$223, (respectively 22$\ldots$223;
122$\ldots$22), is formed in $ H_{j+1}$, if $k+1<m$ and $i\geq 1$,
(respectively $k+1<m$ and $i=0$; $k=m-1$ and $i\geq 1$), with a
2-covered subsequence of $f(H_{j+1})$ exactly on the same
columns of each corresponding 0-covered subsequence of
$f(H_j)$. This bell effect restricts the assignment of $f$ to
the other positions of $ H_{j+1}$. Those of them which are to
receive new label 4 are uniquely determined by steps 1(A)(b2)(i)
and 1(B)(b2)(i), for 2-covered subsequences of $f(H_j)$
immediately at the left of 1-labelled and at
the right of 3-labelled vertices, respectively, so that
corresponding 4-covered subsequences of $f(H_{j+1})$ appear
on the same columns. Step 2 copies each subsequence 123
from positions $(i,j),(i+1,j),(i+2,j)$ to positions
$(i,j+1),(i+1,j+1),(i+2,j+1)$, ($0\leq i<m-2)$. (Here an
alternative option of labels 040 for $(i,j+1),(i+1,j+1),(i+2,j+1)$
would yield two labels 2 at distance 2 on level $ H_{j+2}$, a
contradiction to the definition of a PDS). Now, in each one of the
three decision instances indicated in steps 3 and 4, there are two
possible options: $\alpha$ and $\beta$. The absence of vertices
with label 0 at some level $ H_{j+1}$, signaled by the nullity of
$\tau(j)$ in step 5, yields a PDS in ${\Gamma}_{m,j+2}$ formed
by the vertices that have label 2. Otherwise, steps 1 to 5 must be
repeated with $j:=j+1$. \qfd

\begin{thm}%\fica assim?
If $n$ is finite, then an exhaustive search $\overline{\Theta}$ in
$T_\Theta$ yields all the {\rm PDS}s in ${\mathcal
G}_{m,j}\subseteq{\Gamma}_{m,n}$ with fixed IAVS $S'$ of $ H_0$.
\end{thm}

\proof $\Theta$ may be applied as a sub-procedure of
$\overline{\Theta}$. In each branching of $T_\Theta$,
the exhaustive search
$\overline{\Theta}$ considers first $\alpha$ and then $\beta$. Eventually, step 5 of
$\Theta$ settles $\tau =0$, thus yielding a leaf of
$T_\Theta$. This way, $\overline{\Theta}$ produces all PDSs $S$
of subgraphs ${\Gamma}_{m,j}$ of
${\Gamma}_{m,n}$ with common side $ H_0$ and  such
that $S\cap H_0=S'$. \qfd

\begin{cor}
If $n$ is finite, then $\overline{\Theta}$ %of {\rm Theorem 3.2}
yields
all {\rm PDS}s in ${\Gamma}_{m,n}$ with fixed IAVS $S'$ of $ H_0$.
\end{cor}

\proof The PDSs in the statement are
given by the leaves of $T_\Theta$ at the $n$-th tree row-level,
defined in terms of row order in ${\Gamma}_{m,n}$,
not branching level. \qfd

\begin{prop} Any strategy $\sigma$ of $\Theta$
spends $O(m\times n)$ time to determine whether there exists a {\rm PDS}
$S(m,j+2,S')$ as in
{\rm Theorem 3.1}, even if $j+2=n$.
\end{prop}

\proof No more than five composed complete passes of each $ H_j$
in steps 1 to 5 are performed by $\Theta$. Thus, $5n$ is the number of passes the levels
$H_j$ must be subjected to, in order to determine whether there exists a PDS $S(m.n.S')$
\qfd

\begin{cor}
The running time of $\overline{\Theta}$ in {\rm Theorem 3.2} is
$O(2^{m+n})$.
\end{cor}

\proof An upper bound on the size of the binary decision tree
of $\overline{\Theta}$ is $2^{m+n}$. This and Corollary 3.3 yield
the assertion. \qfd

\section{Questions, greed and threaded pruning}

We restrict to the following case of Question 1.3,
accompanied by two additional questions.

\begin{ques} Given
$1\leq m\in\ZZ$ and an IAVS $S'\subset H_0$, does
there exist $n\in\ZZ$ with $1<n$ such that ${\Gamma}_{m,n}$
contains a {\rm PDS} $S$ with $S\cap H_0=S'$?
\end{ques}

\begin{ques}
In case that $n$
exists, which is the minimal value it attains?\end{ques}

\begin{ques} Which is the spectrum of values of $n$?
\end{ques}

\subsection{The greedy strategy}

Particular strategies $\sigma$
of $\Theta$ may be considered, such as the $\alpha$-strategy.
This is the greedy
strategy for $\Theta$, as in the proof of Theorem 1.4, since by selecting only $\alpha$
when an option is requested,
it proceeds to set locally the least number of labels 0 at any stage of the
running in trying to expose a level $H_j$ of ${\Gamma}_{m,\infty}$
that would not contain dominating vertices of $H_{j+1}$, thus
`greedily pushing' for $n:=j$. A PDS obtained via the $\alpha$-strategy is called a
greedy PDS.

From the proof of Theorem 1.4, it is easy to see that
a modification of the $\alpha$-strategy at its step 5, Subsection 2.4, yields a
periodic PDS $S$ in ${\Gamma}_{m,\infty}$. The resulting modified
algorithm, let us call it $\Theta'$, displays an irreducible period in
${\Gamma}_{m,\infty}$, regardless of the existence of
PDSs in grid graphs ${\Gamma}_{m,n}$. %Solo \alpha-strategy?

In search for an affirmative answer to Question 4.1,
a speedier algorithm in \cite{DD} %http://home.coqui.net/dejterij/sec7.pdf
uses the PDS-arrays that we present in Section 5, below,
aiming at a possible classification of PDSs in grid and band graphs via $\Theta'$.
This was used to show that: {\bf(a)} the lowest $m$
for which no greedy PDS in a finite ${\Gamma}_{m,n}$ exists, for some IAVS
$S'$ of $ H_0$ via $\Theta'$, is 15, and that {\bf(b)} all greedy
PDSs in ${\Gamma}_{m,\infty}$, with $m\le 21$ are embeddable in
toroidal quotients in grid graphs.

It would remain to use a combination of continuing $\alpha$'s and
$\beta$'s to establish particular PDSs in grid graphs ${\mathcal
G}_{m,n}$. This takes us to Subsection 4.2, where a finite pruned
threaded decision tree is obtained that has its leftmost, greedy,
path treated in the link mentioned above.

\subsection{Pruning $T_\Theta$ into a finite threaded tree}

Denote by $\Theta_{S'}$ the algorithm dealing with the fixed IAVS
$S'\subset V(H_0)$. The tree $T_{\Theta_{S'}}$ can be pruned into
a finite threaded tree $T'_{\Theta_{S'}}$ by means of the
following {\it Binary-Tree-Modifying Rule}:

In the running of $\Theta_{S'}(\alpha,\beta)$, as well as in the
corresponding stepwise construction of $T_{\Theta_{S'}}$, each
time a PDS-period $P$ in ${\Gamma}_{m,\infty}$ is engendered,
we opt to stop and avoid the option $\tau\in\{\alpha,\beta\}$ that
leads back to
$P$ in the particular growth stage of $T_{\Theta_{S'}}$. %continuar!
This corresponds to pruning the edge $e$ that descends via option
$\tau$. Instead, a {\it thread} is set from the initial vertex of
$e$ to the initial vertex of $P$. This thread still represents
$\tau$, and its addition to $T_{\Theta_{S'}}$ generates a
well-defined cycle representing $P$.

%We have to ensure that
A continuation of this pruning procedure on $T_{\Theta_{S'}}$ yields a finite threaded tree
$T'_{\Theta_{S'}}$, for any width $m$ of ${\Gamma}_{m,\infty}$ and IAVS
$S'$. To see this, one can think about $T_{\Theta_{S'}}$ as a binary tree with
$\alpha=0$ and $\beta=1$. The lexicographic list $L$ of irreducible descending paths,
(those that are not a concatenation of two or more copies of a shorter path), has the form
(0, 1, 01, 10, 001, 010, 011, 100, 101, 110, 0001, $\ldots$),
and no member path in it is uniformly composed by just two or more 0's or 1's.

For fixed width $m$ of the band graph
${\Gamma}_{m,\infty}$, consider all the possible PDS-slices
$S''=S\cap H_i$ of PDSs in ${\Gamma}_{m,\infty}$.
Any of them can be taken as an initial condition instead of just an IAVS $S'$ for
$\Theta(0,1)$. For every initial $S''$, periodic repetition
of each irreducible descending path will produce through $\Theta'=\Theta'_{S''}$,
a period $P_k(S'')$ formed by the concatenation of a number $k$ of different PDS-slices.
That means that $k$ depends on $S''$, so we may write $k=k(S'')$.

For fixed width $m$, let $\kappa(m)$ be the maximum $k(S'')$, with
$S''$ varying over the (finite)
collection of initial PDS-slices of $\Theta'=\Theta'_m$.
We assert that any infinite descending path $Q$ of a tree $T_{\Theta'_{S''}}$
contains a concatenation of $\kappa(m)$ contiguous copies of $S''$,
for some initial PDS-slice
$S''$ contained in $H_0$. To prove this assertion, we
restrict the infinite lexicographic list $L$ of irreducible
descending paths to the list $L_Q$ of those paths which are contained in $Q$, and
show that $L_Q$ is finite.

Assume that $L=\{p_0,p_1,p_2,\ldots\}$, where $p_0=0$,
$p_1=1$, $p_2=01$, etc. There exists an integer $\mu>0$ such that, for every integer
$x>\mu$, the binary sequence $p_x$ contains a concatenation of $\kappa(m)$ contiguous
copies of some $p_y$, where $\mu\ge y\in\ZZ$. This shows that $L_Q$ is finite. Then,
$Q$ is as required. Now, the union $U$ of the families of irreducible descending paths
contained in all infinite descending paths $Q$ is finite as well.
This implies the following theorem.

\begin{thm}
Let $S'$ be an IAVS contained in $H_0$.
Then the tree $T_{\Theta_{S'}}$ can be pruned into a finite threaded tree
$T'_{\Theta_{S'}}$ by means of the above-mentioned Binary-Tree-Modifying Rule.
\qfd\end{thm}
\begin{cor}
Each finite threaded tree $T_{\Theta_{S'}}$ determines a
well-defined clos\-ed walk $L_{\Theta_{S'}}$ inducing in
${\Gamma}_{m,\infty}$ all possible PDS-slices $S\cap H_i$ obtainable
from $S'$.
\end{cor}
\proof $L_{\Theta_{S'}}$ is obtained by descending $T'_{\Theta_{S'}}$
lexicographically, detouring via every thread found in the covered itinerary of
$T_{\Theta'_{S'}}$,
so every cycle determined by a period of $\Theta_{S'}$ appears exactly once in $L_\Theta$.
In fact, in the process of descending lexicographically $T'_{\Theta'_{S'}}$, each time
the procedure indicates the endvertex $v$ of a thread $e$, the itinerary of $L_{\Theta'_{S'}}$
must completely cover any cycle corresponding to an irreducible descending path in $U$
that ends up with $e$. This may repeat some of the vertices
and edges of $T'_{\Theta'_{S'}}$ in $L_{\Theta'_{S'}}$. For example, the vertex $v$ may be
endvertex of more than one thread $e$. In that case, the procedure should proceed
lexicographically: first with the threads resulting from  the left subtree at $v$, and only
then with the threads resulting from the right subtree at $v$, etc.
This  will account for all possible PDS-slices in $T'_{\Theta'_{S'}}$.
\qfd

\section{PDS-arrays}

For $0<m,n\in\ZZ$, there exists a unique graph
${\Gamma}'_{m,n}$ isomorphic to ${\Gamma}_{m+2,n+2}$ and having
${\Gamma}_{m,n}$ as the subgraph induced by
its interior vertices.
Let the {\it associated graph} ${\Gamma}_{m,n}^S$ of a PDS $S$
in ${\Gamma}_{m,n}$ be the subgraph of ${\Gamma}_{m,n}$ induced
by the complement of $S$ in ${\Gamma}_{m,n}$. Let ${\mathcal H}_{m,n}^S$ be the subgraph
of ${\Gamma}'_{m,n}$ induced by the union of its boundary
cycle with ${\Gamma}_{m,n}^S$. This notation extends to $n=\infty$.

The graph ${\mathcal H}_{m,n}^S$ has chordless
cycles, or holes, delimiting rectangles of areas at least 4,
that we call {\it rooms}, and
maximal connected unions of 4-cycles arranged either
horizontally or vertically into rectangles,
that we call {\it ladders}, a particular case of which is a
4-cycle bordered by rooms. (Formal definitions of rooms and ladders are given in the
proof of Theorem 5.1). The totality of ordered pairs formed by the
horizontal and vertical dimensions, (widths and heights), of the
rectangles spanned by these rooms and ladders can be presented
in an array of integer pairs that we call a {\it PDS-array}.
The two PDSs of Figure 1 yield respectively the following
PDS-arrays:
$$\begin{array}{cccccccccccc}
12 & 42 & 41 & 22 & 21 & 32 & 12 & & & 12 & 22 & 31  \\
22 & 31 & 42 & 21 & 23 & 21 & 22 & & & 22 & 11 & 32  \\
21 & 42 & 21 & 32 & 21 & 23 & 21 & & & 21 & 22 & 21  \\
32 & 31 & 23 & 21 & 32 & 21 & 23 & & & 32 & 11 & 22  \\
31 & 32 & 21 & 23 & 21 & 32 & 21 & & & 31 & 22 & 12  \\
43 & 12 & 32 & 21 & 23 & 21 & 32 & & &     &   & \\
41 & 22 & 12 & 32 & 21 & 22 & 31 & & &     &   &
\end{array}$$
where we write $wh$ for the pair $(w,h)$ formed by the width $w$
and the height $h$ of each rectangle spanned by a room or a
ladder of ${\mathcal H}_{m,n}^S$. Notice that
a pair $(w,h)$ represents a ladder if and only if $min\{w,h\}=1$,
and that pairs representing rooms and ladders are alternate in rows and columns
of the array.

More specifically, given $\delta\in\{0,1\}$ and positive integers $m,n,r,s$,
a PDS-array is defined as a collection
${\mathcal A}={\mathcal A}(m,n,r,s,\delta)=
\{(a_{i,j},b_{i,j}):i\in[0,r),j\in[0,s)\}$ of pairs of positive integers
such that:
\begin{enumerate}\item
If $i+j\equiv\delta$ mod 2 then either $a_{i,j}=1$ or $b_{i,j}=1$,
for $i\in[0,r)$, $j\in[0,s)$;
\item $|a_{i+1,j}-a_{i,j}|\le 2$, for $i\in[0,r-1)$, $j\in[0,s]$;
\item $|b_{i,j+1}-b_{i,j}|\le 2$, for $i\in[0,r)$, $j\in[0,s-1)$;
\item $|\Sigma_{j=0}^ka_{i+1,j}-\Sigma_{j=0}^ka_{i,j}|\le 1$, for
$i\in[0,r-1)$, $k\in[0,s)$;
\item $|\Sigma_{i=0}^kb_{i,j+1}-\Sigma_{i=0}^kb_{i,j}|\le 1$, for
$k\in[0,r)$, $j\in[0,s-1)$;
\item $\Sigma_{i=0}^{r-1}a_{i,j}=m+1$, for $j\in[0,s)$;
\item $\Sigma_{j=0}^{s-1}b_{i,j}=n+1$, for $i\in[0,r)$.
\end{enumerate}
The examples of PDS-arrays given above in relation to
Figure 1 are denoted in fact
${\mathcal A}(16,11,7,7,0)$ and ${\mathcal A}(5,7,3,5,0)$.
If the initial condition in each of the
two given examples were read backwards, the first case would be
${\mathcal A}(16,11,7,7,0)$ again; the second one, ${\mathcal
A}(5,7,3,5,0)$.
\begin{figure}
\unitlength=0.60mm
\special{em:linewidth 0.4pt}
\linethickness{0.4pt}
\begin{picture}(205.30,116.45)
\put(39.06,99.23){\circle{2.00}}
\put(49.06,99.23){\circle*{2.00}}
\put(69.06,89.23){\circle*{2.00}}
\put(39.06,89.23){\circle{2.00}}
\put(49.06,89.23){\circle{2.00}}
\put(59.06,89.23){\circle{2.00}}
\put(39.06,79.23){\circle*{2.00}}
\put(49.06,79.23){\circle{2.00}}
\put(59.06,79.23){\circle{2.00}}
\put(69.06,79.23){\circle{2.00}}
\put(39.06,69.23){\circle{2.00}}
\put(49.06,69.23){\circle{2.00}}
\put(59.06,69.23){\circle*{2.00}}
\put(39.06,102.23){\makebox(0,0)[cc]{$_1$}}
\put(49.06,102.23){\makebox(0,0)[cc]{$_2$}}
\put(59.06,102.23){\makebox(0,0)[cc]{$_3$}}
\put(69.06,102.23){\makebox(0,0)[cc]{$_0$}}
\put(39.06,92.23){\makebox(0,0)[cc]{$_0$}}
\put(49.06,92.23){\makebox(0,0)[cc]{$_4$}}
\put(59.06,92.23){\makebox(0,0)[cc]{$_1$}}
\put(69.06,92.23){\makebox(0,0)[cc]{$_2$}}
\put(39.06,82.23){\makebox(0,0)[cc]{$_2$}}
\put(49.06,82.23){\makebox(0,0)[cc]{$_3$}}
\put(59.06,82.23){\makebox(0,0)[cc]{$_0$}}
\put(69.06,82.23){\makebox(0,0)[cc]{$_4$}}
\put(39.06,72.23){\makebox(0,0)[cc]{$_4$}}
\put(49.06,72.23){\makebox(0,0)[cc]{$_1$}}
\put(59.06,72.23){\makebox(0,0)[cc]{$_2$}}
\put(69.06,72.23){\makebox(0,0)[cc]{$_3$}}
\put(24.06,99.23){\makebox(0,0)[cc]{$^0$}}
\put(24.06,89.23){\makebox(0,0)[cc]{$^1$}}
\put(24.06,79.23){\makebox(0,0)[cc]{$^2$}}
\put(24.06,69.23){\makebox(0,0)[cc]{$^3$}}
\put(39.89,114.23){\makebox(0,0)[cc]{$^0$}}
\put(49.06,114.23){\makebox(0,0)[cc]{$^1$}}
\put(59.06,114.23){\makebox(0,0)[cc]{$^2$}}
\put(69.06,114.23){\makebox(0,0)[cc]{$^3$}}
\put(24.06,109.23){\makebox(0,0)[cc]{$_{j=}$}}
\put(29.06,114.23){\makebox(0,0)[cc]{$_{i=}$}}
\emline{28.06}{110.23}{1}{28.06}{110.23}{2}
\emline{28.58}{110.23}{3}{22.36}{116.45}{4}
\emline{28.25}{110.23}{5}{78.25}{110.23}{6}
\emline{28.25}{110.23}{7}{28.25}{60.23}{8}
\put(59.06,99.23){\circle{2.00}}
\put(69.06,99.23){\circle{2.00}}
\put(69.06,69.23){\circle{2.00}}
\put(39.06,39.23){\circle{2.00}}
\put(49.06,39.23){\circle*{2.00}}
\put(69.06,29.23){\circle*{2.00}}
\put(39.06,29.23){\circle{2.00}}
\put(49.06,29.23){\circle{2.00}}
\put(59.06,29.23){\circle{2.00}}
\put(39.06,19.23){\circle*{2.00}}
\put(49.06,19.23){\circle{2.00}}
\put(59.06,19.23){\circle{2.00}}
\put(39.06,9.23){\circle{2.00}}
\put(49.06,9.23){\circle{2.00}}
\put(39.06,42.23){\makebox(0,0)[cc]{$_1$}}
\put(49.06,42.23){\makebox(0,0)[cc]{$_2$}}
\put(59.06,42.23){\makebox(0,0)[cc]{$_3$}}
\put(69.06,42.23){\makebox(0,0)[cc]{$_0$}}
\put(39.06,32.23){\makebox(0,0)[cc]{$_0$}}
\put(49.06,32.23){\makebox(0,0)[cc]{$_4$}}
\put(59.06,32.23){\makebox(0,0)[cc]{$_1$}}
\put(69.06,32.23){\makebox(0,0)[cc]{$_2$}}
\put(39.06,22.23){\makebox(0,0)[cc]{$_2$}}
\put(49.06,22.23){\makebox(0,0)[cc]{$_3$}}
\put(58.06,22.23){\makebox(0,0)[cc]{$_1$}}
\put(69.06,22.23){\makebox(0,0)[cc]{$_2$}}
\put(39.06,12.23){\makebox(0,0)[cc]{$_2$}}
\put(47.06,11.23){\makebox(0,0)[cc]{$_3$}}
\put(57.06,11.23){\makebox(0,0)[cc]{$_1$}}
\put(69.06,12.23){\makebox(0,0)[cc]{$_2$}}
\put(24.06,39.23){\makebox(0,0)[cc]{$^0$}}
\put(24.06,29.23){\makebox(0,0)[cc]{$^1$}}
\put(24.06,19.23){\makebox(0,0)[cc]{$^2$}}
\put(24.06,9.23){\makebox(0,0)[cc]{$^3$}}
\put(39.89,54.23){\makebox(0,0)[cc]{$^0$}}
\put(49.06,54.23){\makebox(0,0)[cc]{$^1$}}
\put(59.06,54.23){\makebox(0,0)[cc]{$^2$}}
\put(69.06,54.23){\makebox(0,0)[cc]{$^3$}}
\put(24.06,49.23){\makebox(0,0)[cc]{$_{j=}$}}
\put(29.06,54.23){\makebox(0,0)[cc]{$_{i=}$}}
\emline{28.06}{50.23}{9}{28.06}{50.23}{10}
\emline{28.58}{50.23}{11}{22.36}{56.45}{12}
\emline{28.25}{50.23}{13}{78.25}{50.23}{14}
\emline{28.25}{50.23}{15}{28.25}{0.23}{16}
\put(59.06,39.23){\circle{2.00}}
\put(69.06,39.23){\circle{2.00}}
\put(69.06,19.23){\circle*{2.00}}
\put(69.06,19.23){\circle{2.00}}
\put(69.06,9.23){\circle*{2.00}}
\put(59.06,9.23){\circle{2.00}}
\put(39.06,9.23){\circle*{2.00}}
\put(65.43,81.50){\oval(22.08,10.00)[]}
\put(65.43,21.50){\oval(22.08,10.00)[]}
\put(101.25,96.42){\makebox(0,0)[cc]{$_1$}}
\put(101.25,91.42){\makebox(0,0)[cc]{$_0$}}
\put(101.25,86.42){\makebox(0,0)[cc]{$_2$}}
\put(101.25,81.42){\makebox(0,0)[cc]{$_4$}}
\put(106.25,96.42){\makebox(0,0)[cc]{$_2$}}
\put(106.25,91.42){\makebox(0,0)[cc]{$_4$}}
\put(106.25,86.42){\makebox(0,0)[cc]{$_3$}}
\put(106.25,81.42){\makebox(0,0)[cc]{$_1$}}
\put(111.25,96.42){\makebox(0,0)[cc]{$_3$}}
\put(111.25,91.42){\makebox(0,0)[cc]{$_1$}}
\put(111.25,86.42){\makebox(0,0)[cc]{$_0$}}
\put(111.25,81.42){\makebox(0,0)[cc]{$_3$}}
\put(116.25,96.42){\makebox(0,0)[cc]{$_0$}}
\put(116.25,91.42){\makebox(0,0)[cc]{$_2$}}
\put(116.25,86.42){\makebox(0,0)[cc]{$_4$}}
\put(116.25,81.42){\makebox(0,0)[cc]{$_3$}}
\put(132.03,88.42){\makebox(0,0)[cc]{$\Rightarrow$}}
\put(86.03,88.42){\makebox(0,0)[cc]{$\Rightarrow$}}
\emline{146.25}{101.42}{17}{146.25}{91.42}{18}
\emline{141.25}{91.42}{19}{156.25}{91.42}{20}
\emline{156.25}{101.42}{21}{156.25}{86.42}{22}
\emline{151.25}{86.42}{23}{166.25}{86.42}{24}
\emline{151.25}{91.42}{25}{151.25}{76.42}{26}
\emline{141.25}{101.42}{27}{141.25}{76.42}{28}
\emline{141.25}{76.42}{29}{166.25}{76.42}{30}
\emline{166.25}{76.42}{31}{166.25}{101.42}{32}
\emline{166.25}{101.42}{33}{141.25}{101.42}{34}
\emline{141.25}{96.42}{35}{146.25}{96.42}{36}
\emline{156.25}{96.42}{37}{166.25}{96.42}{38}
\emline{161.25}{101.42}{39}{161.25}{96.42}{40}
\emline{161.25}{86.42}{41}{161.25}{76.42}{42}
\emline{161.25}{81.42}{43}{166.25}{81.42}{44}
\emline{141.25}{81.42}{45}{151.25}{81.42}{46}
\emline{146.25}{81.42}{47}{146.25}{76.42}{48}
\put(151.25,96.42){\circle*{2.00}}
\put(146.25,86.42){\circle*{2.00}}
\put(156.25,81.42){\circle*{2.00}}
\put(161.25,91.42){\circle*{2.00}}
\put(101.25,36.42){\makebox(0,0)[cc]{$_1$}}
\put(101.25,31.42){\makebox(0,0)[cc]{$_0$}}
\put(101.25,26.42){\makebox(0,0)[cc]{$_2$}}
\put(101.25,21.42){\makebox(0,0)[cc]{$_2$}}
\put(106.25,36.42){\makebox(0,0)[cc]{$_2$}}
\put(106.25,31.42){\makebox(0,0)[cc]{$_4$}}
\put(106.25,26.42){\makebox(0,0)[cc]{$_3$}}
\put(106.25,21.42){\makebox(0,0)[cc]{$_3$}}
\put(111.25,36.42){\makebox(0,0)[cc]{$_3$}}
\put(111.25,31.42){\makebox(0,0)[cc]{$_1$}}
\put(111.25,26.42){\makebox(0,0)[cc]{$_1$}}
\put(111.25,21.42){\makebox(0,0)[cc]{$_1$}}
\put(116.25,36.42){\makebox(0,0)[cc]{$_0$}}
\put(116.25,31.42){\makebox(0,0)[cc]{$_2$}}
\put(116.25,26.42){\makebox(0,0)[cc]{$_2$}}
\put(116.25,21.42){\makebox(0,0)[cc]{$_2$}}
\put(132.03,28.42){\makebox(0,0)[cc]{$\Rightarrow$}}
\put(86.03,28.42){\makebox(0,0)[cc]{$\Rightarrow$}}
\emline{146.25}{41.42}{49}{146.25}{31.42}{50}
\emline{141.25}{31.42}{51}{156.25}{31.42}{52}
\emline{151.25}{31.42}{53}{151.25}{16.42}{54}
\emline{141.25}{41.42}{55}{141.25}{16.42}{56}
\emline{141.25}{16.42}{57}{166.25}{16.42}{58}
\emline{166.25}{41.42}{59}{141.25}{41.42}{60}
\emline{141.25}{36.42}{61}{146.25}{36.42}{62}
\emline{156.25}{36.42}{63}{166.25}{36.42}{64}
\emline{161.25}{41.42}{65}{161.25}{36.42}{66}
\put(151.25,36.42){\circle*{2.00}}
\put(146.25,26.42){\circle*{2.00}}
\put(161.25,31.42){\circle*{2.00}}
\put(187.20,92.61){\makebox(0,0)[cc]{$_{12}$}}
\put(187.20,87.61){\makebox(0,0)[cc]{$_{22}$}}
\put(187.20,82.61){\makebox(0,0)[cc]{$_{21}$}}
\put(196.01,92.61){\makebox(0,0)[cc]{$_{22}$}}
\put(196.01,87.61){\makebox(0,0)[cc]{$_{11}$}}
\put(196.01,82.61){\makebox(0,0)[cc]{$_{22}$}}
\put(205.30,92.61){\makebox(0,0)[cc]{$_{21}$}}
\put(205.30,87.61){\makebox(0,0)[cc]{$_{22}$}}
\put(205.30,82.61){\makebox(0,0)[cc]{$_{12}$}}
\put(187.20,32.61){\makebox(0,0)[cc]{$_{12}$}}
\put(187.20,27.61){\makebox(0,0)[cc]{$_{23}$}}
\put(196.01,32.61){\makebox(0,0)[cc]{$_{22}$}}
\put(196.01,27.61){\makebox(0,0)[cc]{$_{13}$}}
\put(205.30,32.61){\makebox(0,0)[cc]{$_{21}$}}
\put(205.30,27.61){\makebox(0,0)[cc]{$_{24}$}}
\put(175.84,88.42){\makebox(0,0)[cc]{$\Rightarrow$}}
\put(175.84,28.42){\makebox(0,0)[cc]{$\Rightarrow$}}
\put(146.25,21.42){\circle*{2.00}}
\put(161.25,21.42){\circle*{2.00}}
\put(161.25,26.42){\circle*{2.00}}
\emline{166.25}{41.42}{67}{166.25}{16.42}{68}
\emline{156.25}{41.42}{69}{156.25}{16.42}{70}
\emline{156.25}{21.42}{71}{151.25}{21.42}{72}
\emline{151.25}{26.42}{73}{156.25}{26.42}{74}
\end{picture}
\caption{Two examples of applications of $\Theta$ to PDS-arrays}
\end{figure}
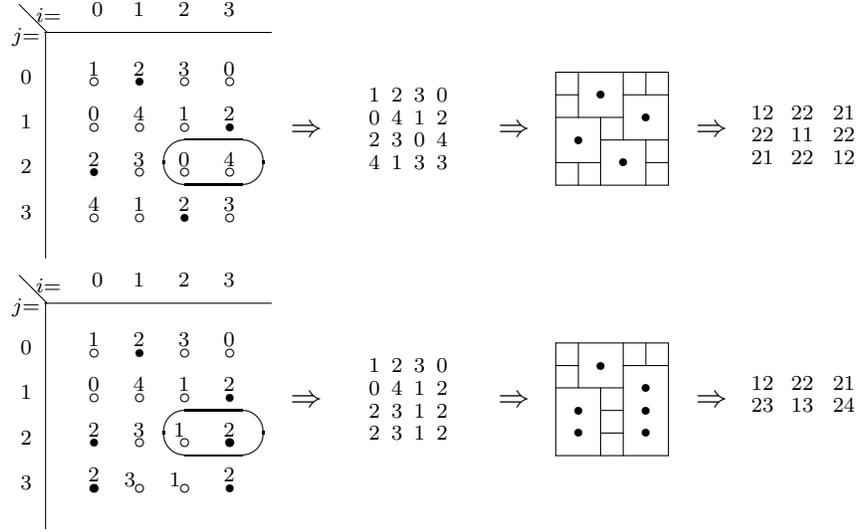
We extend the definition of ${\mathcal A}(m,n,r,s,\delta)$ to
${\mathcal A}(m,\infty,r,\infty,\delta)=$
%\noindent
$\{(a_{i,j},b_{i,j}):i\in[0,r),j\in[0,\infty)\}$, %\subseteq [1,\infty)^2$
ruled solely by items (1-6), with the subindices
$j\in[0,s]$ in items (3-6) replaced by $j\in[0,\infty)$.

Figure 2 contains two elementary examples of the assignment
produced by Theorem 5.1 below, namely ${\mathcal A}(4,4,2,2,0)$ and ${\mathcal
A}(4,4,2,1,0)$, both for the same initial condition $S'$, given by $f(H_0)=1230$,
and having just one differing option, $\alpha,\beta$, respectively.
The figure assigns, to each one of the two runnings of $\Theta$ depicted, as in Figure 1,
on the left side of the figure, and subsequently advancing to the right via the
assignment symbol $\Rightarrow$, the following objects: {\bf(a)}
the corresponding array $f({\Gamma}_{4,4})$, {\bf(b)} the graph ${\mathcal H}_{4,4}^S$,
and finally {\bf(c)} the
PDS-array containing the widths and heights of ${\mathcal H}_{4,4}^S$.

The first (top) example yields the only existing isolated PDS in a grid graph
${\Gamma}_{m,n}$ with $min\{m,n\}>2$
(up to symmetry) \cite{LS}. The second (bottom) example is purposely
continued one more level, once a PDS was obtained with three
levels, which illustrates the following general fact:
\begin{enumerate}\item[] In the running of $\Theta$ with a strategy $\sigma$,
say $S'$ is a PDS-slice in $H_j$
such that $f(H_j)$ has no entries in $\{0,4\}$.
Then, a PDS $S$ with $S'=S\cap H_j$ exists in
${\Gamma}_{m,j+1}$, and if $S''\equiv S'$ is a PDS-slice in $H_{j+1}$,
then a PDS $\overline{S}\supset S$ such that $S''=\overline{S}\cap H_{j+1}$
exists in ${\Gamma}_{m,j+2}$.
\end{enumerate}
If the initial condition in the
two examples of Figure 2 were read backwards, from right to left, so that
labels 1 and 3 were exchanged, the associated PDS-arrays would be
${\mathcal A}(4,4,2,2,0)$ and ${\mathcal A}(4,4,2,1,0)$,
respectively.

\begin{thm} There is an injective assignment from the family of {\rm PDS}s
in grid graphs ${\Gamma}_{m,n}$ into the family of {\rm
PDS}-arrays.\end{thm}

\proof Recall that a function $f$ as in Section 2 can always be defined on a
grid graph ${\Gamma}_{m,n}$ containing a PDS $S$ by assigning labels:

\[ f(i,j) = \left\{ \begin{array}{cccc}
0, & \mbox{if} & (i,j+1) & \in  S;\\
1, & \mbox{if} & (i+1,j) & \in  S;\\
2, & \mbox{if} & (i,j)   & \in  S;\\
3, & \mbox{if} & (i-1,j) & \in  S;\\
4, & \mbox{if} & (i,j-1) & \in  S; \end{array} \right. \]
where $(i,j)\in V({\Gamma}_{m,n})$.
The disposition of such labels for the vertices of the
4-cycles in ${\Gamma}_{m,n}$ are of two different types:
{\bf (a)} those with at least one label 2:
%namely
$$\begin{array}
{lllllllllllll}
00 & 01 & 04 & 12 & 12 & 12 & 22 & 22 & 23 & 23 & 23 & 30 & 40 \\
22 & 23 & 23 & 04 & 12 & 34 & 22 & 44 & 23 & 40 & 41 & 12 & 12,
\end{array}$$
yielding {\it good} 4-cycles, and {\bf (b)}
those having no label 2:
$$\begin{array}
{lllllllllllllll}
31 & 31 & 31 & 31 & 34 & 41 & 41 & 41 & 41 & 44 & 44 & 44 & 44 \\
00 & 01 & 30 & 31 & 30 & 00 & 01 & 30 & 31 & 00 & 01 & 30 & 01,
\end{array}$$
yielding {\it bad} 4-cycles. We extend these notions of
{\it goodness} and {\it badness} to the 4-cycles of
${\Gamma}'_{m,n}$ by attaching to the vertices $v$ of the
boundary path ${\Gamma}'{m,n}\setminus{\Gamma}_{m,n}$
of ${\Gamma}'_{m,n}$ an $f(v)\notin\{0,1,2,3,4\}$, say $f(v)=5$.

For a fixed PDS $S$ in a grid graph ${\Gamma}_{m,n}$: {\bf (a)} the union
of the squares spanned by the good 4-cycles of ${\Gamma}_{m,n}$
splits into maximal rectangles
that we call {\it rooms}; %, whose sides are of length $>1$;
{\bf (b)} the
union of the squares spanned by the bad 4-cycles of $S$ splits into maximal rectangles
that we call {\it ladders}.

If a room or a ladder $R$ has maximal horizontal, (vertical),
common path length $\ell$ in ${\Gamma}_{m,n}^S$, then we say that
the {\it horizontal}, ({\it vertical}), {\it extended length} of $R$ is its width,
(height), as a member of ${\mathcal H}_{m,n}^S$, which equals: {\bf
(a)} $\ell$, if $R$ does not touch a vertical, (horizontal),
boundary path of ${\Gamma}_{m,n}$; {\bf (b)} $\ell +1$, if $R$
touches exactly one vertical, (horizontal), boundary path of
${\Gamma}_{m,n}$; {\bf (c)} $\ell +2$, if $R$ touches both
vertical, (horizontal), boundary paths of ${\Gamma}_{m,n}$.

Every room, (ladder), has exactly one ladder, (room), adjacent to
each side not in the boundary of ${\mathcal
G}'_{m,n}$. The ordered pairs formed by the
horizontal and vertical extended lengths of each room and
ladder can be arranged by means of these room-ladder adjacency instances,
yielding a PDS-array $\phi(S)$. An assignment $\phi$ takes place,
which is as required in the statement. \qfd

Observe that different ladders of a grid graph ${\mathcal H}_{m,n}^S$ are
always disjoint, and that each room $C$ is adjacent to exactly one
room $\ne C$ at each corner of $C$ not in the boundary of
${\Gamma}_{m,n}$ along exactly one edge, either horizontal or vertical.

\section{Application to total perfect codes}

A PDS whose induced components are 1-cubes is called a
{\it total perfect code} (TPC). \cite{KG} shows that
a grid ${\Gamma}_{m,n}$ with
$min\{m,n\}>1$ contains a TPC if and only if
$m\equiv 0$ (mod 2) and $n\equiv -3,-1$ or 1 (mod $m+1$).

An application of $\Theta$ leads to Theorem 6.1 on the
existence of only one TPC in ${\Lambda}$ restricting to TPCs in grid graphs.
We start by generating a TPC in ${\Gamma}_{m,m+2}$ for each even $m>1$.
Let $S'$ a vertex set of even
cardinality. Then, the labeling $f(H_0)$ does not have labels 0; the label 2 appears in
subsequences of the form 1223 in such a way that if the 1's and 3's are interchanged
then backward reading preserves $f(H_0)$.
For increasing values of $m>1$, $f(H_0)$ looks like
$22, 1223, 223122, 12231223, 2231223122, 122312231223, \ldots$
A strategy $\gamma$ useful in this application of $\Theta$
consists of the following selections:
\begin{enumerate}
\item in Subsection 2.3, if $k=1$ and $f(\ell,j-1)$ is 0 or undefined,
then select $\alpha$; otherwise, select $\beta$, where $\ell=0$ in Step 3 and $\ell=i+1$
in Step 4;
\item Check Step 5 of Subsection 2.4 only if $j>0$, with $\tau$ replaced by
$$\tau'=|\{i\in[0,m):f(i,j+1)=f(i,0)\}|.$$
\end{enumerate}
Observe that item 2 yields a modified Step $5'$.
Again, this produces a period $P$ of length $m+1$,
by repeating $f(H_0)$ as $f(H_{m+1})$,
Figure 3 (see upper part) illustrates the array version of a TPC in
${\Gamma}_{m,m+2}$, for $m=2,4,6,8,10$.
But the original Step 5 makes
$\Theta$ stop earlier, yielding, for $m$ even $>2$, a TPC in ${\Gamma}_{m,m}$.

Denote $f=f_m$ and $m'=m+2$, for $m$ even $>1$.
A transformation $\Phi$, indicated in the upper part of Figure 3 by the symbol $\Rightarrow$,
consists in rotating $f_m({\Gamma}_{m,m+2})$ $90\deg$ clockwise
around its center and subjected simultaneously to the
permutation $(0341)$, with the resulting array fitting into
the rectangular box depicted over
$f_{m'}({\Gamma}_{m',m'-2})$, to the right of the symbol $\Rightarrow$, for $m=2,4,6,8$.
Thus, $\Phi$ maps $f_m({\Gamma}_{m,m+2})$ bijectively onto the array enclosed in the
mentioned box. On the other hand,
each depicted array $f_m({\Gamma}_{m,m+2})$ is
translated downwards,
(via $\Downarrow$), into its correspondent graph ${\mathcal H}_{m,m+2}^S$,
shown in the lower part of Figure 3, for $m=2,4,6,8,10$.
This produces the following PDS-arrays:

$$\begin{array}{ccccc}
32 & 12\,\,32\,\,12 & 32\,\,12\,\,32 & 12\,\,32\,\,12\,\,32\,\,12 & 32\,\,12\,\,32\,\,12\,\,32\,\,\\
31 & 23\,\,13\,\,23 & 21\,\,32\,\,21 & 23\,\,12\,\,32\,\,12\,\,23 & 21\,\,32\,\,12\,\,32\,\,21 \\
32 & 12\,\,32\,\,12 & 23\,\,31\,\,23 & 12\,\,23\,\,13\,\,23\,\,21 & 23\,\,21\,\,32\,\,21\,\,23 \\
   &          & 21\,\,32\,\,21 & 23\,\,12\,\,32\,\,12\,\,24 & 21\,\,23\,\,31\,\,23\,\,21 \\
   &          & 32\,\,12\,\,32 & 12\,\,32\,\,12\,\,32\,\,12 & 23\,\,21\,\,32\,\,21\,\,23  \\
   &          &          &                & 21\,\,32\,\,12\,\,32\,\,21 \\
   &          &          &                & 32\,\,12\,\,32\,\,12\,\,32 \\
\end{array}$$

\begin{figure}
\unitlength=0.58mm
\special{em:linewidth 0.4pt}
\linethickness{0.4pt}
\begin{picture}(195.78,130.00)
\emline{140.78}{55.00}{1}{195.78}{55.00}{2}
\emline{145.78}{55.00}{3}{145.78}{50.00}{4}
\emline{140.78}{50.00}{5}{150.78}{50.00}{6}
\put(145.78,45.00){\circle*{2.00}}
\put(145.78,40.00){\circle*{2.00}}
\emline{140.78}{35.00}{7}{150.78}{35.00}{8}
\emline{140.78}{30.00}{9}{150.78}{30.00}{10}
\emline{145.78}{35.00}{11}{145.78}{30.00}{12}
\put(145.78,25.00){\circle*{2.00}}
\put(145.78,20.00){\circle*{2.00}}
\emline{140.78}{15.00}{13}{150.78}{15.00}{14}
\emline{145.78}{15.00}{15}{145.78}{10.00}{16}
\emline{140.78}{10.00}{17}{195.78}{10.00}{18}
\emline{150.78}{55.00}{19}{150.78}{10.00}{20}
\put(145.78,60.00){\circle*{2.00}}
\put(150.78,60.00){\circle*{2.00}}
\emline{155.78}{65.00}{21}{155.78}{55.00}{22}
\emline{160.78}{65.00}{23}{160.78}{55.00}{24}
\emline{155.78}{60.00}{25}{160.78}{60.00}{26}
\put(145.78,5.00){\circle*{2.00}}
\put(150.78,5.00){\circle*{2.00}}
\emline{155.78}{10.00}{27}{155.78}{0.00}{28}
\emline{160.78}{10.00}{29}{160.78}{0.00}{30}
\emline{155.78}{5.00}{31}{160.78}{5.00}{32}
\put(155.78,50.00){\circle*{2.00}}
\put(160.78,50.00){\circle*{2.00}}
\emline{165.78}{55.00}{33}{165.78}{45.00}{34}
\emline{170.78}{55.00}{35}{170.78}{45.00}{36}
\emline{165.78}{50.00}{37}{170.78}{50.00}{38}
\put(175.78,50.00){\circle*{2.00}}
\put(180.78,50.00){\circle*{2.00}}
\emline{185.78}{55.00}{39}{185.78}{10.00}{40}
\put(155.78,15.00){\circle*{2.00}}
\put(160.78,15.00){\circle*{2.00}}
\emline{165.78}{20.00}{41}{165.78}{10.00}{42}
\emline{170.78}{20.00}{43}{170.78}{10.00}{44}
\emline{165.78}{15.00}{45}{170.78}{15.00}{46}
\put(175.78,15.00){\circle*{2.00}}
\put(180.78,15.00){\circle*{2.00}}
\emline{150.78}{45.00}{47}{185.78}{45.00}{48}
\emline{185.78}{45.00}{49}{185.78}{45.00}{50}
\emline{150.78}{20.00}{51}{185.78}{20.00}{52}
\emline{185.78}{20.00}{53}{185.78}{20.00}{54}
\emline{155.78}{45.00}{55}{155.78}{40.00}{56}
\emline{150.78}{40.00}{57}{160.78}{40.00}{58}
\emline{160.78}{40.00}{59}{160.78}{40.00}{60}
\emline{160.78}{40.00}{61}{160.78}{40.00}{62}
\emline{160.78}{45.00}{63}{160.78}{20.00}{64}
\emline{155.78}{20.00}{65}{155.78}{25.00}{66}
\emline{150.78}{25.00}{67}{160.78}{25.00}{68}
\put(155.78,35.00){\circle*{2.00}}
\put(155.78,30.00){\circle*{2.00}}
\emline{175.78}{20.00}{69}{175.78}{45.00}{70}
\emline{175.78}{45.00}{71}{175.78}{45.00}{72}
\emline{175.78}{40.00}{73}{185.78}{40.00}{74}
\emline{180.78}{45.00}{75}{180.78}{40.00}{76}
\emline{175.78}{25.00}{77}{185.78}{25.00}{78}
\emline{180.78}{25.00}{79}{180.78}{20.00}{80}
\put(180.78,35.00){\circle*{2.00}}
\put(180.78,30.00){\circle*{2.00}}
\emline{160.78}{35.00}{81}{175.78}{35.00}{82}
\emline{160.78}{30.00}{83}{175.78}{30.00}{84}
\emline{165.78}{35.00}{85}{165.78}{30.00}{86}
\emline{170.78}{35.00}{87}{170.78}{30.00}{88}
\put(165.78,40.00){\circle*{2.00}}
\put(170.78,40.00){\circle*{2.00}}
\put(165.78,25.00){\circle*{2.00}}
\put(170.78,25.00){\circle*{2.00}}
\emline{190.78}{55.00}{89}{190.78}{50.00}{90}
\emline{185.78}{50.00}{91}{195.78}{50.00}{92}
\put(190.78,45.00){\circle*{2.00}}
\put(190.78,40.00){\circle*{2.00}}
\emline{185.78}{35.00}{93}{195.78}{35.00}{94}
\emline{185.78}{30.00}{95}{195.78}{30.00}{96}
\emline{190.78}{35.00}{97}{190.78}{30.00}{98}
\put(190.78,25.00){\circle*{2.00}}
\put(190.78,20.00){\circle*{2.00}}
\emline{185.78}{15.00}{99}{195.78}{15.00}{100}
\emline{190.78}{15.00}{101}{190.78}{10.00}{102}
\put(165.78,5.00){\circle*{2.00}}
\put(170.78,5.00){\circle*{2.00}}
\emline{175.78}{10.00}{103}{175.78}{0.00}{104}
\emline{180.78}{10.00}{105}{180.78}{0.00}{106}
\emline{175.78}{5.00}{107}{180.78}{5.00}{108}
\put(185.78,5.00){\circle*{2.00}}
\put(190.78,5.00){\circle*{2.00}}
\put(165.78,60.00){\circle*{2.00}}
\put(170.78,60.00){\circle*{2.00}}
\emline{175.78}{65.00}{109}{175.78}{55.00}{110}
\emline{180.78}{65.00}{111}{180.78}{55.00}{112}
\emline{175.78}{60.00}{113}{180.78}{60.00}{114}
\put(185.78,60.00){\circle*{2.00}}
\put(190.78,60.00){\circle*{2.00}}
\emline{140.78}{65.00}{115}{195.78}{65.00}{116}
\emline{140.78}{0.00}{117}{195.78}{0.00}{118}
\emline{140.78}{65.00}{119}{140.78}{0.00}{120}
\emline{195.78}{65.00}{121}{195.78}{0.00}{122}
\put(95.78,45.00){\circle*{2.00}}
\put(95.78,40.00){\circle*{2.00}}
\emline{100.78}{45.00}{123}{105.78}{45.00}{124}
\put(110.78,45.00){\circle*{2.00}}
\put(115.78,45.00){\circle*{2.00}}
\emline{120.78}{45.00}{125}{125.78}{45.00}{126}
\put(130.78,45.00){\circle*{2.00}}
\put(130.78,40.00){\circle*{2.00}}
\emline{100.78}{40.00}{127}{125.78}{40.00}{128}
\emline{90.78}{60.00}{129}{135.78}{60.00}{130}
\emline{110.78}{40.00}{131}{110.78}{25.00}{132}
\emline{115.78}{25.00}{133}{115.78}{40.00}{134}
\emline{110.78}{35.00}{135}{115.78}{35.00}{136}
\emline{110.78}{30.00}{137}{115.78}{30.00}{138}
\emline{100.78}{25.00}{139}{125.78}{25.00}{140}
\put(105.78,35.00){\circle*{2.00}}
\put(105.78,30.00){\circle*{2.00}}
\put(120.78,35.00){\circle*{2.00}}
\put(120.78,30.00){\circle*{2.00}}
\put(110.78,20.00){\circle*{2.00}}
\put(115.78,20.00){\circle*{2.00}}
\emline{105.78}{50.00}{141}{105.78}{40.00}{142}
\emline{120.78}{50.00}{143}{120.78}{40.00}{144}
\emline{105.78}{25.00}{145}{105.78}{15.00}{146}
\emline{120.78}{25.00}{147}{120.78}{15.00}{148}
\emline{110.78}{15.00}{149}{110.78}{5.00}{150}
\emline{115.78}{15.00}{151}{115.78}{5.00}{152}
\emline{110.78}{10.00}{153}{115.78}{10.00}{154}
\emline{90.78}{5.00}{155}{135.78}{5.00}{156}
\emline{90.78}{60.00}{157}{90.78}{5.00}{158}
\emline{135.78}{60.00}{159}{135.78}{5.00}{160}
\put(100.78,10.00){\circle*{2.00}}
\put(105.78,10.00){\circle*{2.00}}
\put(120.78,10.00){\circle*{2.00}}
\put(125.78,10.00){\circle*{2.00}}
\emline{95.78}{15.00}{161}{95.78}{5.00}{162}
\emline{90.78}{10.00}{163}{95.78}{10.00}{164}
\emline{130.78}{15.00}{165}{130.78}{5.00}{166}
\emline{130.78}{10.00}{167}{135.78}{10.00}{168}
\emline{90.78}{15.00}{169}{135.78}{15.00}{170}
\emline{90.78}{50.00}{171}{135.78}{50.00}{172}
\emline{110.78}{60.00}{173}{110.78}{50.00}{174}
\emline{115.78}{60.00}{175}{115.78}{50.00}{176}
\emline{110.78}{55.00}{177}{115.78}{55.00}{178}
\put(100.78,55.00){\circle*{2.00}}
\put(105.78,55.00){\circle*{2.00}}
\put(120.78,55.00){\circle*{2.00}}
\put(125.78,55.00){\circle*{2.00}}
\emline{95.78}{60.00}{179}{95.78}{50.00}{180}
\emline{90.78}{55.00}{181}{95.78}{55.00}{182}
\emline{130.78}{60.00}{183}{130.78}{50.00}{184}
\emline{130.78}{55.00}{185}{135.78}{55.00}{186}
\emline{100.78}{50.00}{187}{100.78}{15.00}{188}
\emline{100.78}{20.00}{189}{105.78}{20.00}{190}
\emline{120.78}{20.00}{191}{125.78}{20.00}{192}
\emline{125.78}{15.00}{193}{125.78}{50.00}{194}
\emline{125.78}{35.00}{195}{135.78}{35.00}{196}
\emline{90.78}{35.00}{197}{100.78}{35.00}{198}
\emline{95.78}{35.00}{199}{95.78}{30.00}{200}
\emline{90.78}{30.00}{201}{100.78}{30.00}{202}
\emline{125.78}{30.00}{203}{135.78}{30.00}{204}
\emline{130.78}{30.00}{205}{130.78}{35.00}{206}
\put(95.78,25.00){\circle*{2.00}}
\put(95.78,20.00){\circle*{2.00}}
\put(130.78,25.00){\circle*{2.00}}
\put(130.78,20.00){\circle*{2.00}}
\emline{85.78}{55.00}{207}{85.78}{10.00}{208}
\emline{50.78}{10.00}{209}{50.78}{55.00}{210}
\emline{75.78}{20.00}{211}{75.78}{45.00}{212}
\emline{60.78}{45.00}{213}{60.78}{20.00}{214}
\emline{50.78}{10.00}{215}{85.78}{10.00}{216}
\emline{85.78}{20.00}{217}{50.78}{20.00}{218}
\emline{50.78}{45.00}{219}{85.78}{45.00}{220}
\emline{85.78}{55.00}{221}{50.78}{55.00}{222}
\put(55.78,50.00){\circle*{2.00}}
\put(60.78,50.00){\circle*{2.00}}
\put(80.78,50.00){\circle*{2.00}}
\put(75.78,50.00){\circle*{2.00}}
\emline{70.78}{45.00}{223}{70.78}{55.00}{224}
\emline{65.78}{55.00}{225}{65.78}{45.00}{226}
\emline{65.78}{50.00}{227}{70.78}{50.00}{228}
\emline{75.78}{40.00}{229}{85.78}{40.00}{230}
\emline{80.78}{40.00}{231}{80.78}{45.00}{232}
\emline{55.78}{45.00}{233}{55.78}{40.00}{234}
\emline{50.78}{40.00}{235}{60.78}{40.00}{236}
\emline{60.78}{25.00}{237}{50.78}{25.00}{238}
\emline{55.78}{25.00}{239}{55.78}{20.00}{240}
\emline{75.78}{25.00}{241}{85.78}{25.00}{242}
\emline{80.78}{25.00}{243}{80.78}{20.00}{244}
\emline{65.78}{20.00}{245}{65.78}{10.00}{246}
\emline{65.78}{15.00}{247}{70.78}{15.00}{248}
\emline{70.78}{20.00}{249}{70.78}{10.00}{250}
\emline{75.78}{30.00}{251}{60.78}{30.00}{252}
\emline{60.78}{35.00}{253}{75.78}{35.00}{254}
\emline{70.78}{35.00}{255}{70.78}{30.00}{256}
\emline{65.78}{30.00}{257}{65.78}{35.00}{258}
\put(55.78,15.00){\circle*{2.00}}
\put(60.78,15.00){\circle*{2.00}}
\put(75.78,15.00){\circle*{2.00}}
\put(80.78,15.00){\circle*{2.00}}
\put(55.78,35.00){\circle*{2.00}}
\put(55.78,30.00){\circle*{2.00}}
\put(80.78,35.00){\circle*{2.00}}
\put(80.78,30.00){\circle*{2.00}}
\put(65.78,40.00){\circle*{2.00}}
\put(70.78,40.00){\circle*{2.00}}
\put(65.78,25.00){\circle*{2.00}}
\put(70.78,25.00){\circle*{2.00}}
\emline{20.78}{45.00}{259}{25.78}{45.00}{260}
\put(30.78,45.00){\circle*{2.00}}
\put(35.78,45.00){\circle*{2.00}}
\emline{40.78}{45.00}{261}{45.78}{45.00}{262}
\emline{20.78}{40.00}{263}{45.78}{40.00}{264}
\emline{30.78}{40.00}{265}{30.78}{25.00}{266}
\emline{35.78}{25.00}{267}{35.78}{40.00}{268}
\emline{30.78}{35.00}{269}{35.78}{35.00}{270}
\emline{30.78}{30.00}{271}{35.78}{30.00}{272}
\emline{20.78}{25.00}{273}{45.78}{25.00}{274}
\put(25.78,35.00){\circle*{2.00}}
\put(25.78,30.00){\circle*{2.00}}
\put(40.78,35.00){\circle*{2.00}}
\put(40.78,30.00){\circle*{2.00}}
\put(30.78,20.00){\circle*{2.00}}
\put(35.78,20.00){\circle*{2.00}}
\emline{25.78}{50.00}{275}{25.78}{40.00}{276}
\emline{40.78}{50.00}{277}{40.78}{40.00}{278}
\emline{25.78}{25.00}{279}{25.78}{15.00}{280}
\emline{40.78}{25.00}{281}{40.78}{15.00}{282}
\emline{20.78}{20.00}{283}{25.78}{20.00}{284}
\emline{40.78}{20.00}{285}{45.78}{20.00}{286}
\emline{45.78}{15.00}{287}{45.78}{50.00}{288}
\emline{45.78}{50.00}{289}{20.78}{50.00}{290}
\emline{20.78}{15.00}{291}{20.78}{50.00}{292}
\emline{45.78}{15.00}{293}{20.78}{15.00}{294}
\put(55.78,115.00){\makebox(0,0)[cc]{$_4$}}
\put(55.78,110.00){\makebox(0,0)[cc]{$_0$}}
\put(55.78,105.00){\makebox(0,0)[cc]{$_2$}}
\put(55.78,100.00){\makebox(0,0)[cc]{$_2$}}
\put(55.78,95.00){\makebox(0,0)[cc]{$_4$}}
\put(55.78,90.00){\makebox(0,0)[cc]{$_0$}}
\put(55.78,85.00){\makebox(0,0)[cc]{$_2$}}
\put(60.78,115.00){\makebox(0,0)[cc]{$_4$}}
\put(60.78,110.00){\makebox(0,0)[cc]{$_1$}}
\put(60.78,105.00){\makebox(0,0)[cc]{$_3$}}
\put(60.78,100.00){\makebox(0,0)[cc]{$_3$}}
\put(60.78,95.00){\makebox(0,0)[cc]{$_1$}}
\put(60.78,90.00){\makebox(0,0)[cc]{$_0$}}
\put(60.78,85.00){\makebox(0,0)[cc]{$_2$}}
\put(65.78,115.00){\makebox(0,0)[cc]{$_0$}}
\put(65.78,110.00){\makebox(0,0)[cc]{$_2$}}
\put(65.78,105.00){\makebox(0,0)[cc]{$_4$}}
\put(65.78,100.00){\makebox(0,0)[cc]{$_0$}}
\put(65.78,95.00){\makebox(0,0)[cc]{$_2$}}
\put(65.78,90.00){\makebox(0,0)[cc]{$_4$}}
\put(65.78,85.00){\makebox(0,0)[cc]{$_3$}}
\put(70.78,115.00){\makebox(0,0)[cc]{$_0$}}
\put(70.78,110.00){\makebox(0,0)[cc]{$_2$}}
\put(70.78,105.00){\makebox(0,0)[cc]{$_4$}}
\put(70.78,100.00){\makebox(0,0)[cc]{$_0$}}
\put(70.78,95.00){\makebox(0,0)[cc]{$_2$}}
\put(70.78,90.00){\makebox(0,0)[cc]{$_4$}}
\put(70.78,85.00){\makebox(0,0)[cc]{$_1$}}
\put(75.78,90.00){\makebox(0,0)[cc]{$_0$}}
\put(75.78,85.00){\makebox(0,0)[cc]{$_2$}}
\put(80.78,90.00){\makebox(0,0)[cc]{$_0$}}
\put(80.78,85.00){\makebox(0,0)[cc]{$_2$}}
\put(145.78,120.00){\makebox(0,0)[cc]{$_0$}}
\put(145.78,115.00){\makebox(0,0)[cc]{$_2$}}
\put(145.78,110.00){\makebox(0,0)[cc]{$_2$}}
\put(145.78,105.00){\makebox(0,0)[cc]{$_4$}}
\put(145.78,100.00){\makebox(0,0)[cc]{$_0$}}
\put(145.78,95.00){\makebox(0,0)[cc]{$_2$}}
\put(145.78,90.00){\makebox(0,0)[cc]{$_2$}}
\put(145.78,85.00){\makebox(0,0)[cc]{$_4$}}
\put(145.78,80.00){\makebox(0,0)[cc]{$_0$}}
\put(145.78,75.00){\makebox(0,0)[cc]{$_2$}}
\put(150.78,120.00){\makebox(0,0)[cc]{$_1$}}
\put(150.78,115.00){\makebox(0,0)[cc]{$_3$}}
\put(150.78,110.00){\makebox(0,0)[cc]{$_3$}}
\put(150.78,105.00){\makebox(0,0)[cc]{$_1$}}
\put(150.78,100.00){\makebox(0,0)[cc]{$_1$}}
\put(150.78,95.00){\makebox(0,0)[cc]{$_3$}}
\put(150.78,90.00){\makebox(0,0)[cc]{$_3$}}
\put(150.78,85.00){\makebox(0,0)[cc]{$_1$}}
\put(150.78,80.00){\makebox(0,0)[cc]{$_0$}}
\put(150.78,75.00){\makebox(0,0)[cc]{$_2$}}
\put(95.78,120.00){\makebox(0,0)[cc]{$_0$}}
\put(95.78,115.00){\makebox(0,0)[cc]{$_2$}}
\put(95.78,110.00){\makebox(0,0)[cc]{$_2$}}
\put(95.78,105.00){\makebox(0,0)[cc]{$_4$}}
\put(95.78,100.00){\makebox(0,0)[cc]{$_0$}}
\put(95.78,95.00){\makebox(0,0)[cc]{$_2$}}
\put(95.78,90.00){\makebox(0,0)[cc]{$_2$}}
\put(95.78,85.00){\makebox(0,0)[cc]{$_4$}}
\put(95.78,80.00){\makebox(0,0)[cc]{$_1$}}
\put(155.78,120.00){\makebox(0,0)[cc]{$_2$}}
\put(155.78,115.00){\makebox(0,0)[cc]{$_4$}}
\put(155.78,110.00){\makebox(0,0)[cc]{$_0$}}
\put(155.78,105.00){\makebox(0,0)[cc]{$_2$}}
\put(155.78,100.00){\makebox(0,0)[cc]{$_2$}}
\put(155.78,95.00){\makebox(0,0)[cc]{$_4$}}
\put(155.78,90.00){\makebox(0,0)[cc]{$_0$}}
\put(155.78,85.00){\makebox(0,0)[cc]{$_2$}}
\put(155.78,80.00){\makebox(0,0)[cc]{$_4$}}
\put(155.78,75.00){\makebox(0,0)[cc]{$_3$}}
\put(100.78,120.00){\makebox(0,0)[cc]{$_4$}}
\put(100.78,115.00){\makebox(0,0)[cc]{$_3$}}
\put(100.78,110.00){\makebox(0,0)[cc]{$_3$}}
\put(100.78,105.00){\makebox(0,0)[cc]{$_1$}}
\put(100.78,100.00){\makebox(0,0)[cc]{$_1$}}
\put(100.78,95.00){\makebox(0,0)[cc]{$_3$}}
\put(100.78,90.00){\makebox(0,0)[cc]{$_3$}}
\put(100.78,85.00){\makebox(0,0)[cc]{$_0$}}
\put(100.78,80.00){\makebox(0,0)[cc]{$_2$}}
\put(160.78,120.00){\makebox(0,0)[cc]{$_2$}}
\put(160.78,115.00){\makebox(0,0)[cc]{$_4$}}
\put(160.78,110.00){\makebox(0,0)[cc]{$_1$}}
\put(160.78,105.00){\makebox(0,0)[cc]{$_3$}}
\put(160.78,100.00){\makebox(0,0)[cc]{$_3$}}
\put(160.78,95.00){\makebox(0,0)[cc]{$_1$}}
\put(160.78,90.00){\makebox(0,0)[cc]{$_0$}}
\put(160.78,85.00){\makebox(0,0)[cc]{$_2$}}
\put(160.78,80.00){\makebox(0,0)[cc]{$_4$}}
\put(160.78,75.00){\makebox(0,0)[cc]{$_1$}}
\put(105.78,120.00){\makebox(0,0)[cc]{$_4$}}
\put(105.78,115.00){\makebox(0,0)[cc]{$_1$}}
\put(105.78,110.00){\makebox(0,0)[cc]{$_0$}}
\put(105.78,105.00){\makebox(0,0)[cc]{$_2$}}
\put(105.78,100.00){\makebox(0,0)[cc]{$_2$}}
\put(105.78,95.00){\makebox(0,0)[cc]{$_4$}}
\put(105.78,90.00){\makebox(0,0)[cc]{$_1$}}
\put(105.78,85.00){\makebox(0,0)[cc]{$_0$}}
\put(105.78,80.00){\makebox(0,0)[cc]{$_2$}}
\put(165.78,120.00){\makebox(0,0)[cc]{$_3$}}
\put(165.78,115.00){\makebox(0,0)[cc]{$_0$}}
\put(165.78,110.00){\makebox(0,0)[cc]{$_2$}}
\put(165.78,105.00){\makebox(0,0)[cc]{$_4$}}
\put(165.78,100.00){\makebox(0,0)[cc]{$_0$}}
\put(165.78,95.00){\makebox(0,0)[cc]{$_2$}}
\put(165.78,90.00){\makebox(0,0)[cc]{$_4$}}
\put(165.78,85.00){\makebox(0,0)[cc]{$_3$}}
\put(165.78,80.00){\makebox(0,0)[cc]{$_0$}}
\put(165.78,75.00){\makebox(0,0)[cc]{$_2$}}
\put(110.78,120.00){\makebox(0,0)[cc]{$_0$}}
\put(110.78,115.00){\makebox(0,0)[cc]{$_2$}}
\put(110.78,110.00){\makebox(0,0)[cc]{$_4$}}
\put(110.78,105.00){\makebox(0,0)[cc]{$_3$}}
\put(110.78,100.00){\makebox(0,0)[cc]{$_3$}}
\put(110.78,95.00){\makebox(0,0)[cc]{$_0$}}
\put(110.78,90.00){\makebox(0,0)[cc]{$_2$}}
\put(110.78,85.00){\makebox(0,0)[cc]{$_4$}}
\put(110.78,80.00){\makebox(0,0)[cc]{$_3$}}
\put(170.78,120.00){\makebox(0,0)[cc]{$_1$}}
\put(170.78,115.00){\makebox(0,0)[cc]{$_0$}}
\put(170.78,110.00){\makebox(0,0)[cc]{$_2$}}
\put(170.78,105.00){\makebox(0,0)[cc]{$_4$}}
\put(170.78,100.00){\makebox(0,0)[cc]{$_0$}}
\put(170.78,95.00){\makebox(0,0)[cc]{$_2$}}
\put(170.78,90.00){\makebox(0,0)[cc]{$_4$}}
\put(170.78,85.00){\makebox(0,0)[cc]{$_1$}}
\put(170.78,80.00){\makebox(0,0)[cc]{$_0$}}
\put(170.78,75.00){\makebox(0,0)[cc]{$_2$}}
\put(115.78,120.00){\makebox(0,0)[cc]{$_0$}}
\put(115.78,115.00){\makebox(0,0)[cc]{$_2$}}
\put(115.78,110.00){\makebox(0,0)[cc]{$_4$}}
\put(115.78,105.00){\makebox(0,0)[cc]{$_1$}}
\put(115.78,100.00){\makebox(0,0)[cc]{$_1$}}
\put(115.78,95.00){\makebox(0,0)[cc]{$_0$}}
\put(115.78,90.00){\makebox(0,0)[cc]{$_2$}}
\put(115.78,85.00){\makebox(0,0)[cc]{$_4$}}
\put(115.78,80.00){\makebox(0,0)[cc]{$_1$}}
\put(175.78,120.00){\makebox(0,0)[cc]{$_2$}}
\put(175.78,115.00){\makebox(0,0)[cc]{$_4$}}
\put(175.78,110.00){\makebox(0,0)[cc]{$_3$}}
\put(175.78,105.00){\makebox(0,0)[cc]{$_1$}}
\put(175.78,100.00){\makebox(0,0)[cc]{$_1$}}
\put(175.78,95.00){\makebox(0,0)[cc]{$_3$}}
\put(175.78,90.00){\makebox(0,0)[cc]{$_0$}}
\put(175.78,85.00){\makebox(0,0)[cc]{$_2$}}
\put(175.78,80.00){\makebox(0,0)[cc]{$_4$}}
\put(175.78,75.00){\makebox(0,0)[cc]{$_3$}}
\put(120.78,120.00){\makebox(0,0)[cc]{$_4$}}
\put(120.78,115.00){\makebox(0,0)[cc]{$_3$}}
\put(120.78,110.00){\makebox(0,0)[cc]{$_0$}}
\put(120.78,105.00){\makebox(0,0)[cc]{$_2$}}
\put(120.78,100.00){\makebox(0,0)[cc]{$_2$}}
\put(120.78,95.00){\makebox(0,0)[cc]{$_4$}}
\put(120.78,90.00){\makebox(0,0)[cc]{$_3$}}
\put(120.78,85.00){\makebox(0,0)[cc]{$_0$}}
\put(120.78,80.00){\makebox(0,0)[cc]{$_2$}}
\put(180.78,120.00){\makebox(0,0)[cc]{$_2$}}
\put(180.78,115.00){\makebox(0,0)[cc]{$_4$}}
\put(180.78,110.00){\makebox(0,0)[cc]{$_0$}}
\put(180.78,105.00){\makebox(0,0)[cc]{$_2$}}
\put(180.78,100.00){\makebox(0,0)[cc]{$_2$}}
\put(180.78,95.00){\makebox(0,0)[cc]{$_4$}}
\put(180.78,90.00){\makebox(0,0)[cc]{$_0$}}
\put(180.78,85.00){\makebox(0,0)[cc]{$_2$}}
\put(180.78,80.00){\makebox(0,0)[cc]{$_4$}}
\put(180.78,75.00){\makebox(0,0)[cc]{$_1$}}
\put(125.78,120.00){\makebox(0,0)[cc]{$_4$}}
\put(125.78,115.00){\makebox(0,0)[cc]{$_1$}}
\put(125.78,110.00){\makebox(0,0)[cc]{$_1$}}
\put(125.78,105.00){\makebox(0,0)[cc]{$_3$}}
\put(125.78,100.00){\makebox(0,0)[cc]{$_3$}}
\put(125.78,95.00){\makebox(0,0)[cc]{$_1$}}
\put(125.78,90.00){\makebox(0,0)[cc]{$_1$}}
\put(125.78,85.00){\makebox(0,0)[cc]{$_0$}}
\put(125.78,80.00){\makebox(0,0)[cc]{$_2$}}
\put(185.78,120.00){\makebox(0,0)[cc]{$_3$}}
\put(185.78,115.00){\makebox(0,0)[cc]{$_1$}}
\put(185.78,110.00){\makebox(0,0)[cc]{$_1$}}
\put(185.78,105.00){\makebox(0,0)[cc]{$_3$}}
\put(185.78,100.00){\makebox(0,0)[cc]{$_3$}}
\put(185.78,95.00){\makebox(0,0)[cc]{$_1$}}
\put(185.78,90.00){\makebox(0,0)[cc]{$_1$}}
\put(185.78,85.00){\makebox(0,0)[cc]{$_3$}}
\put(185.78,80.00){\makebox(0,0)[cc]{$_0$}}
\put(185.78,75.00){\makebox(0,0)[cc]{$_2$}}
\put(130.78,120.00){\makebox(0,0)[cc]{$_0$}}
\put(130.78,115.00){\makebox(0,0)[cc]{$_2$}}
\put(130.78,110.00){\makebox(0,0)[cc]{$_2$}}
\put(130.78,105.00){\makebox(0,0)[cc]{$_4$}}
\put(130.78,100.00){\makebox(0,0)[cc]{$_0$}}
\put(130.78,95.00){\makebox(0,0)[cc]{$_2$}}
\put(130.78,90.00){\makebox(0,0)[cc]{$_2$}}
\put(130.78,85.00){\makebox(0,0)[cc]{$_4$}}
\put(130.78,80.00){\makebox(0,0)[cc]{$_3$}}
\put(190.78,120.00){\makebox(0,0)[cc]{$_0$}}
\put(190.78,115.00){\makebox(0,0)[cc]{$_2$}}
\put(190.78,110.00){\makebox(0,0)[cc]{$_2$}}
\put(190.78,105.00){\makebox(0,0)[cc]{$_4$}}
\put(190.78,100.00){\makebox(0,0)[cc]{$_0$}}
\put(190.78,95.00){\makebox(0,0)[cc]{$_2$}}
\put(190.78,90.00){\makebox(0,0)[cc]{$_2$}}
\put(190.78,85.00){\makebox(0,0)[cc]{$_4$}}
\put(190.78,80.00){\makebox(0,0)[cc]{$_0$}}
\put(190.78,75.00){\makebox(0,0)[cc]{$_2$}}
\put(55.78,120.00){\makebox(0,0)[cc]{$_2$}}
\put(60.78,120.00){\makebox(0,0)[cc]{$_2$}}
\put(65.78,120.00){\makebox(0,0)[cc]{$_3$}}
\put(70.78,120.00){\makebox(0,0)[cc]{$_1$}}
\put(95.78,125.00){\makebox(0,0)[cc]{$_1$}}
\put(100.78,125.00){\makebox(0,0)[cc]{$_2$}}
\put(105.78,125.00){\makebox(0,0)[cc]{$_2$}}
\put(110.78,125.00){\makebox(0,0)[cc]{$_3$}}
\put(115.78,125.00){\makebox(0,0)[cc]{$_1$}}
\put(120.78,125.00){\makebox(0,0)[cc]{$_2$}}
\put(125.78,125.00){\makebox(0,0)[cc]{$_2$}}
\put(130.78,125.00){\makebox(0,0)[cc]{$_3$}}
\put(25.78,110.00){\makebox(0,0)[cc]{$_0$}}
\put(25.78,105.00){\makebox(0,0)[cc]{$_2$}}
\put(25.78,100.00){\makebox(0,0)[cc]{$_2$}}
\put(25.78,95.00){\makebox(0,0)[cc]{$_4$}}
\put(25.78,90.00){\makebox(0,0)[cc]{$_1$}}
\put(30.78,110.00){\makebox(0,0)[cc]{$_4$}}
\put(30.78,105.00){\makebox(0,0)[cc]{$_3$}}
\put(30.78,100.00){\makebox(0,0)[cc]{$_3$}}
\put(30.78,95.00){\makebox(0,0)[cc]{$_0$}}
\put(30.78,90.00){\makebox(0,0)[cc]{$_2$}}
\put(35.78,110.00){\makebox(0,0)[cc]{$_4$}}
\put(35.78,105.00){\makebox(0,0)[cc]{$_1$}}
\put(35.78,100.00){\makebox(0,0)[cc]{$_1$}}
\put(35.78,95.00){\makebox(0,0)[cc]{$_0$}}
\put(35.78,90.00){\makebox(0,0)[cc]{$_2$}}
\put(40.78,110.00){\makebox(0,0)[cc]{$_0$}}
\put(40.78,105.00){\makebox(0,0)[cc]{$_2$}}
\put(40.78,100.00){\makebox(0,0)[cc]{$_2$}}
\put(40.78,95.00){\makebox(0,0)[cc]{$_4$}}
\put(40.78,90.00){\makebox(0,0)[cc]{$_3$}}
\put(25.78,115.00){\makebox(0,0)[cc]{$_1$}}
\put(30.78,115.00){\makebox(0,0)[cc]{$_2$}}
\put(35.78,115.00){\makebox(0,0)[cc]{$_2$}}
\put(40.78,115.00){\makebox(0,0)[cc]{$_3$}}
\put(75.78,115.00){\makebox(0,0)[cc]{$_4$}}
\put(75.78,110.00){\makebox(0,0)[cc]{$_3$}}
\put(75.78,105.00){\makebox(0,0)[cc]{$_1$}}
\put(75.78,100.00){\makebox(0,0)[cc]{$_1$}}
\put(75.78,95.00){\makebox(0,0)[cc]{$_3$}}
\put(80.78,115.00){\makebox(0,0)[cc]{$_4$}}
\put(80.78,110.00){\makebox(0,0)[cc]{$_0$}}
\put(80.78,105.00){\makebox(0,0)[cc]{$_2$}}
\put(80.78,100.00){\makebox(0,0)[cc]{$_2$}}
\put(80.78,95.00){\makebox(0,0)[cc]{$_4$}}
\put(75.78,120.00){\makebox(0,0)[cc]{$_2$}}
\put(80.78,120.00){\makebox(0,0)[cc]{$_2$}}
\put(145.78,130.00){\makebox(0,0)[cc]{$_2$}}
\put(145.78,125.00){\makebox(0,0)[cc]{$_4$}}
\put(150.78,130.00){\makebox(0,0)[cc]{$_2$}}
\put(150.78,125.00){\makebox(0,0)[cc]{$_4$}}
\put(155.78,130.00){\makebox(0,0)[cc]{$_3$}}
\put(155.78,125.00){\makebox(0,0)[cc]{$_0$}}
\put(160.78,130.00){\makebox(0,0)[cc]{$_1$}}
\put(160.78,125.00){\makebox(0,0)[cc]{$_0$}}
\put(165.78,130.00){\makebox(0,0)[cc]{$_2$}}
\put(165.78,125.00){\makebox(0,0)[cc]{$_4$}}
\put(170.78,130.00){\makebox(0,0)[cc]{$_2$}}
\put(170.78,125.00){\makebox(0,0)[cc]{$_4$}}
\put(175.78,130.00){\makebox(0,0)[cc]{$_3$}}
\put(175.78,125.00){\makebox(0,0)[cc]{$_0$}}
\put(180.78,130.00){\makebox(0,0)[cc]{$_1$}}
\put(180.78,125.00){\makebox(0,0)[cc]{$_0$}}
\put(185.78,130.00){\makebox(0,0)[cc]{$_2$}}
\put(185.78,125.00){\makebox(0,0)[cc]{$_4$}}
\put(190.78,130.00){\makebox(0,0)[cc]{$_2$}}
\put(190.78,125.00){\makebox(0,0)[cc]{$_4$}}
\emline{53.78}{93.00}{295}{82.78}{93.00}{296}
\emline{82.78}{93.00}{297}{82.78}{112.00}{298}
\emline{82.78}{112.00}{299}{53.78}{112.00}{300}
\emline{53.78}{112.00}{301}{53.78}{93.00}{302}
\emline{132.78}{117.00}{303}{93.78}{117.00}{304}
\emline{23.78}{107.00}{305}{42.78}{107.00}{306}
\emline{42.78}{107.00}{307}{42.78}{98.00}{308}
\emline{42.78}{98.00}{309}{23.78}{98.00}{310}
\emline{23.78}{98.00}{311}{23.78}{107.00}{312}
\emline{93.78}{117.00}{313}{93.78}{88.00}{314}
\emline{93.78}{88.00}{315}{132.78}{88.00}{316}
\emline{132.78}{88.00}{317}{132.78}{117.00}{318}
\emline{143.78}{83.00}{319}{192.78}{83.00}{320}
\emline{192.78}{83.00}{321}{192.78}{122.00}{322}
\emline{192.78}{122.00}{323}{143.78}{122.00}{324}
\emline{143.78}{122.00}{325}{143.78}{83.00}{326}
\put(32.78,70.00){\makebox(0,0)[cc]{$\Downarrow$}}
\put(67.78,70.00){\makebox(0,0)[cc]{$\Downarrow$}}
\put(112.78,70.00){\makebox(0,0)[cc]{$\Downarrow$}}
\put(167.78,70.00){\makebox(0,0)[cc]{$\Downarrow$}}
\put(5.78,110.00){\makebox(0,0)[cc]{$_2$}}
\put(5.78,105.00){\makebox(0,0)[cc]{$_4$}}
\put(5.78,100.00){\makebox(0,0)[cc]{$_0$}}
\put(5.78,95.00){\makebox(0,0)[cc]{$_2$}}
\put(10.78,110.00){\makebox(0,0)[cc]{$_2$}}
\put(10.78,105.00){\makebox(0,0)[cc]{$_4$}}
\put(10.78,100.00){\makebox(0,0)[cc]{$_0$}}
\put(10.78,95.00){\makebox(0,0)[cc]{$_2$}}
\put(7.78,70.00){\makebox(0,0)[cc]{$\Downarrow$}}
\emline{0.78}{45.00}{327}{0.78}{20.00}{328}
\emline{0.78}{20.00}{329}{15.78}{20.00}{330}
\emline{15.78}{20.00}{331}{15.78}{45.00}{332}
\emline{15.78}{45.00}{333}{0.78}{45.00}{334}
\emline{15.78}{35.00}{335}{0.78}{35.00}{336}
\emline{0.78}{30.00}{337}{15.78}{30.00}{338}
\emline{10.78}{30.00}{339}{10.78}{35.00}{340}
\emline{5.78}{35.00}{341}{5.78}{30.00}{342}
\put(5.78,40.00){\circle*{2.00}}
\put(10.78,40.00){\circle*{2.00}}
\put(10.78,25.00){\circle*{2.00}}
\put(5.78,25.00){\circle*{2.00}}
\put(18.78,103.00){\makebox(0,0)[cc]{$\Rightarrow$}}
\put(48.78,103.00){\makebox(0,0)[cc]{$\Rightarrow$}}
\put(88.78,103.00){\makebox(0,0)[cc]{$\Rightarrow$}}
\put(138.78,103.00){\makebox(0,0)[cc]{$\Rightarrow$}}
\end{picture}
\caption{Arrays $f({\Gamma}_{m,m+2})$ for TPCs and corresponding ${\mathcal H}_{m,m+2}^S$}
\end{figure}
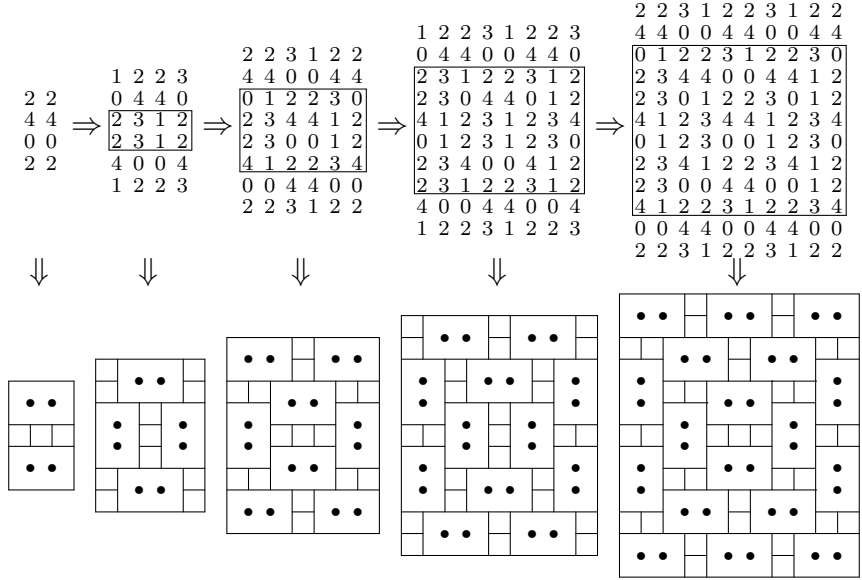

As for the remaining TPCs characterized in \cite{KG},
level $H_2$ filled according to $\gamma$ yields another $S'$ for $\Theta$, by
taking $H_2$ as the initial level $H_0$ via the change of coordinates
$(i,j)\rightarrow(i,j-2)$. This yields a TPC in ${\Gamma}_{m,m-2}$
whose corresponding array $f({\Gamma}_{m,m-2})$ is shown inside
the rectangular box in the upper part of Figure 3, for
$m=2,4,6,8,10$, which forms a subarray of $f({\Gamma}_{m,m+2})$.
By adding to $f({\Gamma}_{m,m-2})$ the first two or last two
rows of $f({\Gamma}_{m,m+2})$, we get two different arrays of
corresponding TPCs in ${\Gamma}_{m,m}$, of which one was cited
above. These two arrays are related by the alphabet permutation
$(04)(13)$, corresponding to a rotation of $180\deg$ around the
center of ${\Gamma}_{m,m}$. However, for $6\le m\equiv 2$ (mod 4), there
is still another TPC in ${\Gamma}_{m,m}$, whose corresponding
$f({\Gamma}_{m,m})$ is obtained by prefixing, (postfixing), the pairs
23, 23, 41, 01, 23, 23, etc., (12, 12, 34, 30, 12, 12, etc.), to the
first, second, $\ldots, m$-th rows of $f({\Gamma}_{m-2,m})$.
Thus, the $S'$ for this TPC is non symmetric, but
$\gamma$ and $\Theta''$ still apply. All these TPCs may be continued
by $\gamma$ and $\Theta''$ with periodicity of minimum period length
$m+1$.

By changing to centered coordinates,
the grid graphs ${\Gamma}_{m,m+2}$, for all even $m>0$,
may be considered each one, say ${\Gamma}_{m,m+2}$,
positioned inside the next one, ${\Gamma}_{m+2,m+4}$,
in a concentric manner in ${\Lambda}$
in such a way that the union of all resulting concentrically positioned TPCs in ${\Gamma}_{m,m+2}$ in ${\Lambda}$ constitute a TPC $S_1$
of ${\Lambda}$. This is clear from Figure 4, because by rotating $90\deg$ each depicted
${\Gamma}_{m,m+2}$ and then translating the resulting grid graph
concentrically inside the next grid graph, that is ${\Gamma}_{m+2,m+4}$, the rotated TPC
obtained from ${\Gamma}_{m,m+2}$ is a concentric
part of the TPC in ${\Gamma}_{m+2,m+4}$. Thus, in the limit for increasing $m$, the
claimed TPC $S_1$ in ${\Lambda}$ is obtained. Moreover, the complement ${\Lambda}-S_1$
yields a tiling formed by rooms and ladders with respective PDS-array entry
sets (23,32) and (12,21,13,31),
(where precisely one central ladder PDS-array entry in (13,31) exists).
This tiling is aperiodic, as is the Penrose tiling
of \cite{Penrose}, so it has no translational symmetry, but
nevertheless has 4-fold rotational symmetry about vertex $(0,0)$
and double reflective symmetry about the coordinate axes $x=0$ and $y=0$.
These together conform the central symmetries produced by the group $D_8$ (of
symmetries of the central square $[-\frac{1}{2},\frac{1}{2}]\times[-\frac{1}{2},\frac{1}{2}]$
of the ladder with PDS-array entry 13 or 31),
which is in fact the automorphism group of the plane tiling in question.

\begin{thm}
The TPC $S_1$ is the only existing TPC in ${\Lambda}$ which restricts to TPCs in
$m\times n$ grid graphs ${\Gamma}_{m,n}$, where $m$ and $n$ are integers $>2$.
Moreover, the complement ${\Lambda}-S_1$ yields an aperiodic tiling of the plane
(like the Penrose tiling)
whose automorphism group coincides with the group $D_8$ of symmetries of the square
$[-\frac{1}{2},\frac{1}{2}]\times[-\frac{1}{2},\frac{1}{2}]$.
\qfd
\end{thm}

This result is in contrast with the uncountability of TPCs shown in Theorem 1 of \cite{DD1}.
Notice however that TPCs can be obtained in band graphs with finite width $m$
out of the TPCs presented in relation to Figure 3, and in their unions having two
parallel maximal lateral paths at distance $m$ with all its vertices of degree 2.

\vspace*{3mm}

%\noindent{\bf Acknowledgement:} The first author is grateful to Jan Kratochvil, Charles
%University, Prague, for his warm
%hospitality and for having suggested a relation of the TPC $S_1$ in ${\Lambda}$ arising
%from our Figure 3 with the Penrose tiling.

 \end{document}